\documentclass[12pt,reqno]{amsart}

\makeatletter
\def\@settitle{%
	\begin{center}%
		\normalfont\Large\bfseries % <-- Ici tu choisis la typo
		\@title
	\end{center}%
}
\makeatother

\usepackage{amsmath}
\usepackage{mathrsfs}
\usepackage{enumitem}   
\usepackage{graphicx}
\usepackage[colorlinks = true,
linkcolor = blue,
urlcolor  = blue,
citecolor = blue,
anchorcolor = blue]{hyperref}

\newtheorem{theorem}{Theorem}
\newtheorem{definition}{Definition}
\newtheorem{proposition}%[theorem]
{Proposition}
\newtheorem{corollary}%[theorem]
{Corollary}
\newtheorem{remark}{Remark}

\newtheorem{lemma}{Lemma}

\newtheorem{example}{Example}

\newfont{\bb}{msbm10 at 12pt}
\def\pf{{\textit {Proof.} }}

\def\<{\langle}     
\def\>{\rangle}     
\def\div{{\rm div}}

\def\Id{\operatorname{Id}}
\def\End{\operatorname{End}}

\newcommand{\bal}{\begin{align}}      \newcommand{\eal}{\end{align}}
\newcommand{\ba}{\begin{array}}      \newcommand{\ea}{\end{array}}
\newcommand{\bc}{\begin{center}}     \newcommand{\ec}{\end{center}}
\newcommand{\be}{\begin{enumerate}}  \newcommand{\ee}{\end{enumerate}}
\newcommand{\beq}{\begin{eqnarray}}  \newcommand{\eeq}{\end{eqnarray}}
\newcommand{\beQ}{\begin{eqnarray*}} \newcommand{\eeQ}{\end{eqnarray*}}
\newcommand{\bi}{\begin{itemize}}    \newcommand{\ei}{\end{itemize}}
\newcommand{\bt}{\begin{tabular}}    \newcommand{\et}{\end{tabular}}
\newcommand{\bdm}{\begin{displaymath}} \newcommand{\edm}{\end{displaymath}}

\newcommand{\D}{D\!\!\!\!/\,}

\newcommand{\RSB}{\mathbf{S}\!\!\!\!/\,}

\newcommand{\nb}{\nabla\!\!\!\!/\,}
\newcommand{\mult}{c\!\!\!/}

\newcommand{\Ss}{R\!\!\!\!/\,}

\def\qed{\hfill{q.e.d.}\smallskip\smallskip}

\begin{document}
	
	\title{Charged parallel spinors and applications to mass--charge inequalities}

	\author{Simon Raulot}
	\address[Simon Raulot]{Univ Rouen Normandie, CNRS, Normandie Univ, LMRS UMR 6085, F-76000 Rouen, France}
	\email{simon.raulot@univ-rouen.fr}
	
	\begin{abstract}
		We investigate the equality case of the spin positive mass theorem with charge in the Riemannian setting. This leads naturally to the notion of charged parallel spinor, which plays a central role in the analysis of extremal charged manifolds. As an application, we characterize the equality case of the mass--charge inequality in terms of the extremal Reissner--Nordstr\"om geometry for asymptotically flat manifolds with connected boundary and for manifolds with a single asymptotically cylindrical end.
	\end{abstract}
	
	\keywords{}

	%\subjclass{53C20, 52C24, 53C27, 83C57}
	
	\thanks{}
	
	\date{\today}   
	
	\maketitle 
	\tableofcontents
	\pagenumbering{arabic}	
	
	%%%%%%%%%%%%%%%%%%%%%%%%%%%%%%%%%%%%%%%%%%%%%%%%%%%%%%%%%%%%%%%%%%%%%%%%%%%%
	
	\section{Introduction}
	
	%%%%%%%%%%%%%%%%%%%%%%%%%%%%%%%%%%%%%%%%%%%%%%%%%%%%%%%%%%%%%%%%%%%%%%%%%%%%
	
	The positive energy theorem is a fundamental result in mathematical relativity. It asserts that, under suitable asymptotic assumptions and the dominant energy condition, the ADM energy--momentum vector $(E,P)\in\mathbb R^{3,1}$ of an initial data set is future-directed and causal, that is, $E \ge |P|$, with equality if and only if the data arise from Minkowski spacetime. This result was first established by Schoen and Yau using minimal surface techniques~\cite{SchoenYau1,SchoenYau6}, and later reproved by Witten via a spinorial argument~\cite{Witten1}.
	
	Witten’s spinorial method extends to higher dimensions under a spin assumption and offers a robust framework for incorporating additional matter fields. In the presence of an electromagnetic field, this approach leads to mass--charge inequalities within the Einstein--Maxwell theory. Early derivations were obtained by Gibbons and Hull~\cite{GibbonsHull} in a supergravity context, and later by Gibbons, Hawking, Horowitz, and Perry~\cite{GibbonsHawkingHorowitzPerry} in the study of black hole spacetimes, both relying on spinorial techniques. The analytical foundations of these arguments were subsequently clarified and made rigorous by Bartnik and Chru\'sciel~\cite{BartnikChrusciel}.
	
	The analysis of the equality case is considerably more subtle and is closely related to the existence of super-covariantly constant spinors, which impose 
	strong geometric constraints on the underlying manifold. Local classifications of such solutions were obtained by Tod~\cite{Tod83}. At the global level, Chru\'sciel, Reall, and Tod~\cite{ChruscielRealTod} showed that, under suitable assumptions, initial data admitting a super-covariantly constant spinor must belong to the Majumdar--Papapetrou class. These manifolds describe electrostatic configurations of extremal charged black holes in equilibrium and provide the natural geometric models for the equality case of the mass--charge inequality.
	
	In higher dimensions, the spinorial approach remains a natural framework for the study of mass--charge inequalities. In the purely electric setting, the author established in previous work~\cite{Raulot16} the corresponding inequality for charged initial data sets using a modified Witten argument. 
	
	Motivated by the time-symmetric Einstein--Maxwell constraints, we consider triples $(M^n,g,E)$ consisting of a Riemannian manifold $(M^n,g)$ together with a vector field $E\in\Gamma(TM)$ interpreted as an electric field. Such a triple will be referred to as a charged manifold. The dominant energy condition is expressed by the scalar inequality
	\begin{equation}\label{DEC}
		\mu:=\frac{1}{2}\left(R - (n-1)(n-2)|E|^2 - 2(n-1)|\div(E)|\right)\geq 0,
	\end{equation}
	where $R$ denotes the scalar curvature of $(M^n,g)$. In this context, the positive energy theorem reduces to the inequality $m \ge |Q|$ relating the ADM mass and the total charge of an asymptotically flat end. Several geometric approaches have been developed to establish this estimate in the time-symmetric setting. 
	
	In dimension three, Jaracz~\cite{Jaracz} proved such an inequality for initial data with either an outermost minimal boundary or asymptotically cylindrical 
	ends, by adapting the inverse mean curvature flow. A different perspective, based on spacetime harmonic functions, was introduced by Bray, Hirsch, Kazaras, Khuri, and Zhang~\cite{BrayHirschKazarasKhuri}. Their approach provides a unified framework for mass inequalities and, importantly, yields a precise description of the equality case under natural geometric assumptions, showing that extremality is strongly constrains by Majumdar--Papapetrou configurations.
	
	Very recently, McCormick~\cite{McCormick} established the Riemannian mass--charge inequality in arbitrary dimension without the spin assumption as a consequence of a positive mass theorem for the \(X\)-ADM mass. His approach is based on a conformal reduction to the classical positive mass theorem and provides a unified proof of several positive mass type inequalities.
	
	Building on the framework developed in~\cite{Raulot16}, the present paper investigates charged manifolds realizing equality in the corresponding mass--charge inequality, which we shall refer to as \emph{extremal charged manifolds}. A central role is played by charged parallel spinors, which arise naturally in this equality case. Rather than being merely a technical consequence of the spinorial proof, these spinors encode the underlying structure of extremal charged manifolds. Their existence determines the electric field through the spinor norm, gives rise to an electrostatic Einstein--Maxwell system, and is closely related, through a conformal change of metric, to the existence of ordinary parallel spinors. 
	
	In the source-free case, a non-trivial total charge requires the presence of an inner geometric structure allowing for non-vanishing electric flux. Such a structure may arise either through a compact inner boundary or through an asymptotically cylindrical end corresponding to a degenerate horizon. These two geometric settings form the framework of the present work. 
	
	Our approach is guided by the viewpoint that charged parallel spinors provide the natural spinorial counterpart of parallel spinors in the Einstein--Maxwell setting. The rigidity statements proved below show how the existence of charged parallel spinors governs the geometry of extremal charged manifolds under suitable global assumptions.
	
	We first consider asymptotically flat manifolds with compact boundary. This setting already captures the main features of the rigidity mechanism and provides a natural context for introducing the boundary contributions arising in the spinorial method.
	
	Our first main result is the following mass--charge inequality. We refer to Section~\ref{MainDefinitions} for the definitions and notations.
	\begin{theorem}\label{ChargedPMT-Connected boundary}
		Let $n \geq 3$ and let $(M^n,g,E)$ be a complete charged spin manifold containing at least one asymptotically flat end and with a compact and connected boundary $\partial M$ with positive Yamabe invariant. Assume that the dominant energy condition~\eqref{DEC} is satisfied and that
		\begin{equation}\label{YamabeMC}
			H+(n-1)|g(E,\nu)|\leq {\rm Vol}(\partial M,h)^{-\frac{1}{n-1}}
			\sqrt{\frac{n-1}{n-2}\mathcal{Y}(\partial M,[h])}.
		\end{equation}
		Then the ADM mass satisfies $m\geq |Q|$. If $\div(E)=0$, equality occurs if and only if $(M^n,g,E)$ is isometric to the exterior region of a coordinate sphere in an extremal Reissner--Nordstr\"om manifold.
	\end{theorem}
	
	Here $h:=g_{|\partial M}$, $H$ denotes the mean curvature of $\partial M$ with respect to the unit normal pointing toward infinity and $\mathcal{Y}(\partial M,[h])$ is the Yamabe invariant of the induced metric $h$. In dimension $n=3$, the Yamabe invariant of the boundary is a topological quantity, and the previous statement simplifies as follows.
	\begin{theorem}\label{ChargedPMT-Connected-3dim}
		Let $(M^3,g,E)$ be a complete charged oriented manifold containing at least one asymptotically flat end and with an inner boundary $\partial M$ homeomorphic to a $2$-sphere. Assume that the dominant energy condition~\eqref{DEC} is satisfied and that
		\[
		H+2|g(E,\nu)|\leq 4\sqrt{\frac{\pi}{{\rm Vol}(\partial M,h)}}.
		\]
		Then the ADM mass satisfies $m\geq |Q|$. If $\div(E)=0$, equality holds if and only if $(M^3,g,E)$ is isometric to the exterior region of a coordinate sphere in an extremal Reissner--Nordstr\"om manifold.
	\end{theorem}
	
	The proof relies on a spectral boundary condition involving the Dirac operator, combined with lower bounds for its first eigenvalue in terms of the Yamabe invariant. The connectedness assumption is only required for the analysis of the equality case, whereas the mass--charge inequality itself remains valid for manifolds with disconnected boundary.
	
	This motivates the study of manifolds with asymptotically cylindrical ends, which can be viewed as a limiting configuration of disconnected boundary components. Our second main result addresses this more delicate situation.
	\begin{theorem}\label{ChargedPMT-Cylindrical}
		Let $n \geq 3$ and let $(M^n,g,E)$ be a complete charged spin manifold containing one asymptotically flat end and a finite number of asymptotically cylindrical ends. Assume that the dominant energy condition~\eqref{DEC} is satisfied. Then the ADM mass $m$ satisfies $m \geq |Q|$. Moreover, if $\div(E)=0$ and $M$ has exactly one asymptotically cylindrical end, then equality holds if and only if $(M^n,g,E)$ is isometric to an extremal Reissner--Nordstr\"om manifold.
	\end{theorem}
	
	The restriction to a single asymptotically cylindrical end is used only in the analysis of the equality case. 
	
	The proofs of the rigidity statements rely on a common conformal strategy based on charged parallel spinors. In both settings, the existence of such a spinor gives rise, through a conformal change of the initial metric, to a Ricci-flat asymptotically flat manifold. In the compact boundary setting, the rigidity statement then follows directly from the rigidity part of Herzlich's positive mass theorem~\cite{Herzlich1,Herzlich2}. The asymptotically cylindrical case requires a substantially more delicate analysis. There, the main difficulty is to show that the conformal metric compactifies the cylindrical end to an isolated conical singularity, allowing us to apply the rigidity part of the positive mass theorem with isolated conical singularities of Dai--Sun--Wang~\cite{DaiSunWang1,DaiSunWang2}.
	
	The paper is organized as follows. In Section~\ref{MainDefinitions}, we introduce the class of charged manifolds under consideration and recall the notions of asymptotically flat and asymptotically cylindrical ends. Section~\ref{SG-CM} develops the spinorial framework and establishes a modified Schr\"odinger--Lichnerowicz formula adapted to the Einstein--Maxwell setting. Section~\ref{Section-CPS} introduces charged parallel spinors, which play a central role in the analysis of the equality case. We investigate their geometric consequences and establish a conformal characterization relating them to metrics carrying ordinary parallel spinors. In Section~\ref{MP-Example}, we study the Majumdar--Papapetrou family, which provides the natural model geometries associated with this spinorial structure. The mass--charge inequality and its rigidity statement are then established in two different geometric settings. Section~\ref{PMT-CB} treats manifolds with compact boundary, while Section~\ref{PMT-ACE} deals with manifolds having asymptotically cylindrical ends.

	%%%%%%%%%%%%%%%%%%%%%%%%%%%%%%%%%%%%%%%%%%%%%%%%%%%%%%%%%%%%%%%%%%%%%%%%%%%%%%%%%%%%%%%%%%%%%%%%%%%%%%%
	
	\section{Notations and definitions}\label{MainDefinitions}
	
	%%%%%%%%%%%%%%%%%%%%%%%%%%%%%%%%%%%%%%%%%%%%%%%%%%%%%%%%%%%%%%%%%%%%%%%%%%%%%%%%%%%%%%%%%%%%%%%%%%%%%%%
	
	In this section, we introduce the class of charged manifolds considered throughout the paper together with the associated geometric quantities.
	\begin{definition}
		A \emph{charged manifold} is a triple $(M^n,g,E)$ where $M$ is a connected, smooth manifold of dimension $n\ge 3$, $g$ is a Riemannian metric on $M$, and $E \in \Gamma(TM)$ is a vector field on $M$.
	\end{definition}
	
	\noindent
	The vector field $E$ will be referred to as the \emph{electric field}. For simplicity, we assume throughout that $g$ and $E$ are smooth, although the results remain valid under weaker regularity assumptions (see \cite{BartnikChrusciel}). 
	
	The notions of asymptotically flat and asymptotically cylindrical ends play a central role in the later sections of the paper. We therefore recall them here.
	
	\begin{definition}
		A subset \(M_{\mathrm{ext}}\subset M\) is called an \emph{asymptotically flat end} if there exists a diffeomorphism
		\[
		M_{\mathrm{ext}}
		\longrightarrow
		\mathbb R^n\setminus\overline B_1(0)
		\]
		such that
		\begin{equation}\label{AFDecay}
			\left\{
			\begin{aligned}
				|(\nabla^\delta)^j(g-\delta)|_\delta
				&=
				O(|x|^{-\tau-j}),
				\qquad j=0,1,2,\\
				|(\nabla^\delta)^jE|_\delta
				&=
				O(|x|^{-\tau-1-j}),
				\qquad j=0,1.
			\end{aligned}
			\right.
		\end{equation}
		for some \(\tau>\frac{n-2}{2}\). Here \(\nabla^\delta\) denotes the Levi-Civita connection of the Euclidean metric \(\delta\) as well as its natural extensions to tensor fields. We also assume that \(\mu\) and \(\div(E)\) are integrable on \(M\).
	\end{definition}
	
	The \emph{ADM mass} of an asymptotically flat end is defined by
	\[
	m
	=
	\frac{1}{2(n-1)\omega_{n-1}}
	\lim_{r\to\infty}
	\int_{S_r}
	\Big(
	\div_\delta(g-\delta)
	-
	d\operatorname{tr}_\delta(g-\delta)
	\Big)(\bar\nu_r)\,
	d\bar\sigma_r,
	\]
	where \(S_r\) denotes the Euclidean sphere of radius \(r>0\), \(\bar\nu_r\) its outward unit normal, \(d\bar\sigma_r\) the induced Euclidean volume element, and \(\omega_{n-1}\) the volume of the unit \((n-1)\)-sphere. Here \(\div_\delta\) and \(\operatorname{tr}_\delta\) are computed with respect to the Euclidean metric \(\delta\). Although this definition appears to depend on the choice of asymptotically flat coordinates, it defines a geometric invariant  \cite{Bartnik1,Chrusciel1}. 
	In the presence of an electric field \(E\in\Gamma(TM)\), the \emph{total charge} is defined by
	\[
		Q
		=
		\frac{1}{\omega_{n-1}}
		\lim_{r\to\infty}
		\int_{S_r}
		\delta(E,\bar\nu_r)\,
		d\bar\sigma_r .
	\]
	
	We now introduce the notion of an asymptotically cylindrical end, which will play a central role in the analysis of extremal charged manifolds.
	
	\begin{definition}\label{Def-ACE}
		A subset $\mathcal E \subset M$ is called an \emph{asymptotically cylindrical end} if there exist a real number $T_0$, a diffeomorphism
		\[
		\mathcal E \simeq (T_0,\infty)\times \Sigma,
		\]
		where $\Sigma$ is a closed manifold of dimension $n-1$, a Riemannian metric $h$ on $\Sigma$, and a $t$-independent vector field $E_\infty\in\Gamma(\mathcal{E})$ such that
		\begin{equation}\label{ACDecay}
			\left\{
			\begin{aligned}
				|(\nabla^\infty)^j(g-g_\infty)|_{g_\infty}
				&=
				O(e^{-\alpha t}),
				\qquad j=0,1,2,\\
				|(\nabla^\infty)^j(E-E_\infty)|_{g_\infty}
				&=
				O(e^{-\alpha t}),
				\qquad j=0,1,
			\end{aligned}
			\right.
		\end{equation}
		for some $\alpha>0$. Here $g_\infty=dt^2+h$ is the Riemannian product metric on the cylinder and $\nabla^\infty$ denotes its Levi-Civita connection as well as its extensions to tensor fields. 
		
		The manifold $(\Sigma^{n-1},h)$ is called the \emph{limiting cross-section} of the asymptotically cylindrical end.
	\end{definition}
	
	We emphasize that no geometric or topological assumption is imposed on the limiting cross-section $(\Sigma^{n-1},h)$.
	
	%%%%%%%%%%%%%%%%%%%%%%%%%%%%%%%%%%%%%%%%%%%%%%%%%%%%%%%%%%%%%%%%%%%%%%%%%%%%
	
	\section{Spin geometry of charged manifolds}\label{SG-CM}
	
	%%%%%%%%%%%%%%%%%%%%%%%%%%%%%%%%%%%%%%%%%%%%%%%%%%%%%%%%%%%%%%%%%%%%%%%%%%%%
	
	Throughout this section, $(M^n,g,E)$ denotes a charged manifold for which $M$ is spin. Such a triple will be called a \emph{charged spin manifold}.
	
	%%%%%%%%%%%%%%%%%%%%%%%%%%%%%%%%%%%%%%%%%%%%%%%%%%%%%%%%%%%%%%%%%%%%%%%%%%%%
	
	\subsection{Basic facts on spinors}\label{Spinors}
	
	%%%%%%%%%%%%%%%%%%%%%%%%%%%%%%%%%%%%%%%%%%%%%%%%%%%%%%%%%%%%%%%%%%%%%%%%%%%%
	
	We briefly recall the basic notions of spin geometry needed in the sequel and fix our notation. For a more detailed presentation, we refer the reader to \cite{LawsonMichelsohn,BourguignonHijaziMilhoratMoroianu,Ginoux}.
	
	Let $(M^n,g)$ be an $n$-dimensional Riemannian spin manifold. Then there exists a smooth Hermitian vector bundle over $M$, the spinor bundle, denoted by $S_g$, whose sections are called spinor fields. The Hermitian scalar product is denoted by $\<\,,\,\>$. Moreover, the tangent bundle $TM$ acts on $S_g$ by Clifford multiplication $X\otimes \psi\mapsto c(X)\psi$ satisfying
	\[
	c(X)c(Y)+c(Y)c(X)=-2g(X,Y)\,{\rm Id}
	\]
	and this action is skew-Hermitian with respect to the Hermitian metric. The Riemannian Levi-Civita connection $\nabla$ lifts to the spin Levi-Civita connection, still denoted by $\nabla$, which defines a metric covariant derivative on $S_g$ preserving Clifford multiplication. The Dirac operator is the first order elliptic differential operator acting on $S_g$ locally defined by
	\[
	D\varphi=\sum_{j=1}^nc(e_j)\nabla_{e_j}\varphi
	\]
	for $\varphi\in\Gamma(S_g)$. Here and throughout the paper, $(e_1,\dots,e_n)$ denotes a local orthonormal frame on $(M^n,g)$. It satisfies the Schr\"odinger-Lichnerowicz formula 
	\begin{equation}\label{SL-Formula}
		D^2\varphi=\nabla^*\nabla\varphi+\frac{R}{4}\varphi
	\end{equation}
	for all $\varphi\in\Gamma(S_g)$ and where $\nabla^*$ denotes the $L^2$-formal adjoint of $\nabla$. If $\Omega\subset M$ is a bounded domain with smooth boundary, then 
	\begin{equation}\label{Dirac-IPP}
		\int_\Omega\<D\varphi,\psi\>d\mu=\int_\Omega\<\varphi,D\psi\>d\mu+\oint_{\partial\Omega}\<c(\nu)\varphi,\psi\>d\sigma
	\end{equation}
	where $\nu$ denotes the outward unit normal to $\partial\Omega$, and $d\mu$ (resp. $d\sigma$) is the Riemannian volume form of $\Omega$ (resp. $\partial\Omega$).
	
	Now let $\Sigma^{n-1}\subset M$ be an oriented hypersurface with unit normal vector field $\nu$ and induced metric $h$. The hypersurface $\Sigma$ inherits a spin structure and we denote by $\RSB_h:=S_{g|N}$ and by $S_h$ the \emph{extrinsic} and the \emph{intrinsic} spinor bundles over $(\Sigma^{n-1},h)$. These bundles are related by the identification
	\begin{equation}\label{SpinorBundleRestriction}
		\RSB_h \simeq
		\begin{cases}
			S_h & \text{if } n \text{ is even},\\
			S_h \oplus S_h & \text{if } n \text{ is odd}.
		\end{cases}
	\end{equation}
	On $\RSB_h$, one defines the \emph{extrinsic Clifford multiplication} and the \emph{extrinsic spin connection} by
	\begin{equation}\label{CliffordLCI-Extrinsic}
		\begin{cases}
			\mult_h(X)\psi &:=  c(X)c(\nu)\psi, \\
			\nb^h_X \psi &:= \nabla_X \psi + \frac{1}{2} \mult_h(AX)\psi, 
		\end{cases}
	\end{equation}
	for all $X\in\Gamma(T\Sigma)$ and $\psi\in\Gamma(\RSB_h)$, where $A:=\nabla\nu$ denotes the shape operator of $\Sigma$. The corresponding Dirac operator is then as usual locally given by
	\begin{equation}\label{ExtrinsicDirac}
		\D_h\varphi=\sum_{i=1}^{n-1}\mult_h(e_i)\nb^h_{e_i}\varphi
	\end{equation}
	for all $\varphi\in\Gamma(\RSB_h)$. These structures satisfy the usual compatibility relations and, under the identification~\eqref{SpinorBundleRestriction}, one has:
	\begin{equation}\label{CliffordLCI-Identifications}
		\nb^h\simeq
		\begin{cases}
			\nabla^h & \text{if } n \text{ is even},\\
			\nabla^h \oplus \nabla^h & \text{if } n \text{ is odd},
		\end{cases}
		\qquad
		\mult_h \simeq
		\begin{cases}
			c_h & \text{if } n \text{ is even},\\
			c_h \oplus -c_h & \text{if } n \text{ is odd}
		\end{cases}
	\end{equation}
	where $\nabla^h$ and $c_h$ denote respectively the spin Levi-Civita connection and the Clifford multiplication on the intrinsic spinor bundle $S_h$. 
	
	We conclude this subsection by recalling the identification of spinor bundles under a change of metric. Let $\overline{g}$ be another Riemannian metric on $M$. Although the spin structure is fixed, the associated spinor bundle depends on the metric. Following Bourguignon--Gauduchon \cite{BourguignonGauduchon}, there exists a unique positive endomorphism $\beta : TM \to TM$ which is symmetric with respect to $g$ and characterized by $g(X,Y)=\overline{g}(\beta X,\beta Y)$ for all $X,Y\in \Gamma(TM)$. It induces a fiberwise isometric bundle isomorphism (still denoted by $\beta$)
	\[
	\beta : S_g\longrightarrow S_{\overline g},
	\]
	which allows us to identify spinor fields defined with respect to $g$ and $\overline{g}$. When $\overline{g} = e^{2f}g\in[g]$ belongs to the conformal class $[g]$ of $g$, the identification becomes explicit (see, e.g., \cite{BourguignonHijaziMilhoratMoroianu,Ginoux}). In this case, if \(\overline{\nabla}\) denotes the spin Levi-Civita connection associated with \(\overline g\), the spin Levi-Civita connections are related by
	\begin{equation}\label{ConformalLCC}
		\overline{\nabla}_X \beta\psi = \beta\big(\nabla_X\psi-\frac{1}{2}c(X)c(\nabla f)\psi-\frac{1}{2}g(\nabla f,X)\psi\big)
	\end{equation}
	for $X\in\Gamma(TM)$ and $\psi\in\Gamma(S_g)$. As we will see, conformal changes of metric arise naturally in the study of charged parallel spinors.
	
	%%%%%%%%%%%%%%%%%%%%%%%%%%%%%%%%%%%%%%%%%%%%%%%%%%%%%%%%%%%%%%%%%%%%%%%%%%%%
	
	\subsection{The modified connections}
	
	%%%%%%%%%%%%%%%%%%%%%%%%%%%%%%%%%%%%%%%%%%%%%%%%%%%%%%%%%%%%%%%%%%%%%%%%%%%%
	
	We introduce the modified connection
	\[
		\nabla^-_X\varphi
		:=
		\nabla_X\varphi
		+
		\frac{1}{2}c(X)c(E)\varphi
		+
		\frac{n-1}{2}g(E,X)\varphi
	\]
	for all $X\in \Gamma(TM)$ and $\varphi\in\Gamma(S_g)$. The associated Dirac operator is locally given by
	\[
	D^-\varphi
	:=
	\sum_{j=1}^nc(e_j)\nabla^-_{e_j}\varphi.
	\]
	
	The first properties of $D^-$ are summarized in the following lemma.
	\begin{lemma}\label{MD-Properties}
		The operator $D^-$ is a first-order elliptic differential operator satisfying
		\begin{equation}\label{RelationsDOD}
			D^-\varphi
			=
			D\varphi
			-
			\frac{1}{2}c(E)\varphi
		\end{equation}
		and its formal $L^2$-adjoint is
		\begin{equation}\label{FormalAdjoint}
			(D^-)^\ast\varphi
			=
			D\varphi
			+
			\frac{1}{2}c(E)\varphi
		\end{equation}
		for all $\varphi\in\Gamma(S_g)$.
	\end{lemma}
	
	\pf
	The fact that $D^-$ is a first-order elliptic differential operator follows directly from (\ref{RelationsDOD}) since $D^-$ differs from the Dirac operator $D$ by a zeroth-order term. The formulas (\ref{RelationsDOD}) and (\ref{FormalAdjoint}) follow by direct computations using the Clifford rule and the fact that Clifford multiplication by $E$ is pointwise skew-Hermitian. 
	\qed
	
	\begin{remark}\label{IPP-DM}
		If $\Omega$ is a bounded domain with smooth boundary $\partial\Omega$, it follows from~\eqref{Dirac-IPP} that
		\[
		\int_\Omega\<D^-\varphi,\psi\>d\mu
		=
		\int_\Omega\<\varphi,(D^-)^\ast\psi\>d\mu
		+
		\oint_{\partial\Omega}\<c(\nu)\varphi,\psi\>d\sigma
		\]
		for all $\varphi, \psi \in \Gamma(S_{g|\Omega})$. 
	\end{remark}
	
	To derive the Weitzenböck formula for the modified Dirac operator $D^-$, we first compute the $L^2$-formal adjoint of $\nabla^-$.
	\begin{lemma}\label{MC-Properties}
		The $L^2$-formal adjoint of $\nabla^-$ satisfies
		\[
		(\nabla^-)^\ast\nabla^-\varphi  
		=  
		\nabla^\ast\nabla\varphi
		-
		\frac{1}{2}\big[D,c(E)\big]\varphi
		+
		\frac{1}{4}\Big(\big((n-1)(n-2)+1\big)|E|^2+2(n-1)\div(E)\Big)\varphi 
		\]
		for all $\varphi\in\Gamma(S_g)$.  
	\end{lemma}
	
	\pf 
	For $\varphi,\psi\in\Gamma(S_g)$  we compute
	\[
	\<\nabla^-\varphi,\nabla^-\psi\> 
	=  
	\div(\xi)+(1)+(2)+(3) 
	\]
	where $\xi\in\Gamma(TM)$ is the vector field defined by $g(\xi,X)=\<\nabla^-_X\varphi,\psi\>$ for all $X\in\Gamma(TM)$, and
	\[
	\begin{aligned}
		(1) &=
		\sum_{j=1}^n
		\bigl\langle
		-\nabla_{e_j}\nabla^-_{e_j}\varphi,
		\psi
		\bigr\rangle, \\[4pt]
		(2) &=
		\frac{1}{2}
		\sum_{j=1}^n
		\bigl\langle
		\nabla^-_{e_j}\varphi,
		c(e_j)c(E)\psi
		\bigr\rangle, \\[4pt]
		(3) &=
		\frac{n-1}{2}
		\sum_{j=1}^n
		g(E,e_j)
		\bigl\langle
		\nabla^-_{e_j}\varphi,
		\psi
		\bigr\rangle.
	\end{aligned}
	\]
	A direct computation yields
	\[
	(1)  
	=  
	\<\nabla^*\nabla\varphi-\frac{1}{2}D(c(E)\varphi)-\frac{n-1}{2}\div(E)\varphi-\frac{n-1}{2}\nabla_E\varphi,\psi\>.
	\]
	Similarly,
	\[
	(2) 
	=  
	\<\frac{1}{2}c(E)D\varphi+\frac{1}{4}|E|^2\varphi,\psi\>.
	\]
	and 
	\[
	(3) 
	=  
	\<\frac{n-1}{2}\nabla_E\varphi+\frac{(n-1)(n-2)}{4}|E|^2\varphi,\psi\>
	\]
	which yields the desired formula.
	\qed
	
	\begin{remark}
		In the preceding proof, we established the following pointwise identity:
		\[
		\langle (\nabla^-)^*\nabla^-\varphi,\psi\rangle
		=
		\langle \nabla^-\varphi,\nabla^-\psi\rangle
		-
		\div(\xi),
		\]
		where $\xi\in\Gamma(TM)$ is the vector field defined by $g(\xi,X)=\<\nabla^-_X\varphi,\psi\>$ for all $X\in\Gamma(TM)$ and $\varphi,\psi\in\Gamma(S_g)$. If $\Omega$ is a bounded domain with smooth boundary in a complete charged spin manifold $(M^n,g,E)$, then applying the divergence theorem yields
		\[
		\int_\Omega
		\langle (\nabla^-)^\ast\nabla^-\varphi,\psi\rangle
		\, d\mu
		=
		\int_\Omega
		\langle \nabla^-\varphi,\nabla^-\psi\rangle
		\, d\mu
		-
		\oint_{\partial\Omega}
		\langle \nabla^-_\nu\varphi,\psi\rangle
		\, d\sigma.
		\]
	\end{remark}
	
	We now state the main formula of this section. The resulting identity is the Riemannian counterpart of the formula obtained in \cite{Raulot16}.
	\begin{theorem}\label{SLforEM}
		Let $(M^n,g,E)$ be a charged spin manifold. Then
		\[
		(D^-)^\ast D^-\varphi
		=
		(\nabla^-)^*\nabla^-\varphi
		+
		\mathcal{R}^-\varphi
		\]
		for all $\varphi\in\Gamma(S_g)$, where $\mathcal{R}^-$ is the function defined by
		\[
		\mathcal{R}^-
		=
		\frac{1}{4}
		\Big(
		R
		-
		(n-1)(n-2)|E|^2
		+
		2(n-1)\div(E)
		\Big).
		\]
	\end{theorem}
	
	\pf 
	A straightforward computation using Lemma~\ref{MD-Properties} shows that
	\[
	(D^-)^\ast D^-\varphi
	=
	D^2\varphi
	-\frac{1}{2}[D,c(E)]\varphi
	+\frac{1}{4}|E|^2\varphi
	\]
	which, combined with the Schr\"odinger-Lichnerowicz formula~\eqref{SL-Formula}, yields
	\[
	(D^-)^\ast D^-\varphi
	=
	\nabla^*\nabla\varphi
	-\frac{1}{2}[D,c(E)]\varphi
	+\frac{1}{4}(R+|E|^2)\varphi.
	\]
	Applying Lemma~\ref{MC-Properties} to this identity yields the result.
	\qed
	
	Integrating the identity of Theorem~\ref{SLforEM} over a bounded domain yields the following formula.
	\begin{corollary}\label{IntegralVersion}
		Let $\Omega$ be a bounded domain with smooth boundary in a complete charged spin manifold $(M^n,g,E)$. Then for all $\varphi\in\Gamma(S_{g|\Omega})$,
		\[
		\int_\Omega
		\Big(
		|\nabla^-\varphi|^2
		+ \mathcal{R}^-|\varphi|^2
		-|D^-\varphi|^2
		\Big)
		\, d\mu
		=
		\oint_{\partial\Omega}
		\langle L^-_{\nu}\varphi,\varphi\rangle
		\, d\sigma.
		\]
		where $\nu$ is the outward unit normal to $\partial\Omega$ and $L^-_\nu\varphi:=\nabla^-_{\nu}\varphi+c(\nu)D^-\varphi$.
	\end{corollary}
	
	Replacing the electric field $E$ by $-E$ leads to a second modified connection
	\[
	\nabla^+_X\varphi
	:=
	\nabla_X\varphi
	-
	\frac{1}{2}c(X)c(E)\varphi
	-
	\frac{n-1}{2}g(E,X)\varphi
	\]
	for all $X\in\Gamma(TM)$ and $\varphi\in\Gamma(S_g)$. The corresponding Dirac operator then satisfies $D^+=(D^-)^\ast$. Moreover,
	\[
	(D^+)^\ast D^+\varphi
	=
	(\nabla^+)^*\nabla^+\varphi
	+
	\mathcal{R}^+\varphi
	\]
	for all $\varphi\in\Gamma(S_g)$, where $\mathcal{R}^+$ is the function defined by
	\[
	\mathcal{R}^+
	=
	\frac{1}{4}
	\Big(
	R
	-
	(n-1)(n-2)|E|^2
	-
	2(n-1)\div(E)
	\Big).
	\]
	Arguing as above, one obtains the analogue of Corollary~\ref{IntegralVersion}:
	\[
	\int_\Omega
	\Big(
	|\nabla^+\varphi|^2
	+ \mathcal{R}^+|\varphi|^2
	-|D^+\varphi|^2
	\Big)
	\, d\mu
	=
	\oint_{\partial\Omega}
	\langle L^+_{\nu}\varphi,\varphi\rangle
	\, d\sigma.
	\]
	where $L^+_\nu\varphi:=\nabla^+_{\nu}\varphi+c(\nu)D^+\varphi$.
	
	%%%%%%%%%%%%%%%%%%%%%%%%%%%%%%%%%%%%%%%%%%%%%%%%%%%%%%%%%%
	\section{Charged parallel spinors}\label{Section-CPS}
	%%%%%%%%%%%%%%%%%%%%%%%%%%%%%%%%%%%%%%%%%%%%%%%%%%%%%%%%%%
	
	In this section, we consider an arbitrary charged spin manifold $(M^n,g,E)$, without any additional assumptions on its geometry or topology. We introduce the notion of charged parallel spinor, which naturally arises in the analysis of the equality case of the charged Schrödinger--Lichnerowicz formula. As we shall see, the existence of such spinors imposes strong geometric constraints and reveals a close connection with electrostatic Einstein--Maxwell systems and conformal deformations of metrics carrying parallel spinors.
	\begin{definition}
		Let $(M^n,g,E)$ be a charged spin manifold.
		\begin{itemize}
			\item A spinor field $\Phi\in\Gamma(S_g)$ is called a \emph{negative charged parallel spinor} if it is parallel with respect to the connection $\nabla^-$, that is,
			\[
			\nabla^-_X\Phi=0
			\]
			for all $X\in\Gamma(TM)$. Equivalently,
			\begin{equation}\label{NegativeCPS}
				\nabla_X\Phi
				=
				-\frac12 c(X)c(E)\Phi
				-\frac{n-1}{2}g(E,X)\Phi
			\end{equation}
			for all $X\in\Gamma(TM)$.
			
			\item  Similarly, a spinor field $\Phi\in\Gamma(S_g)$ is called a \emph{positive charged parallel spinor} if it is parallel with respect to the connection $\nabla^+$, that is,
			\[
			\nabla^+_X\Phi=0
			\]
			for all $X\in\Gamma(TM)$. Equivalently,
			\begin{equation}\label{PositiveCPS}
				\nabla_X\Phi
				=
				\frac{1}{2}c(X)c(E)\Phi
				+
				\frac{n-1}{2}g(E,X)\Phi
			\end{equation}
			for all $X\in\Gamma(TM)$.
			
			\item A spinor field is called a \emph{charged parallel spinor} if it is either a positive or a negative charged parallel spinor.
		\end{itemize}
	\end{definition}
	
	\begin{example}
		If $E=0$, charged parallel spinors reduce to classical parallel spinors. It is well known (see \cite{WangMcK1}) that, in the simply connected case, this includes the Euclidean space $\mathbb{R}^n$, Calabi--Yau manifolds, hyper-Kähler manifolds, $7$-dimensional manifolds with holonomy $\mathrm{G}_2$, $8$-dimensional manifolds with holonomy $\mathrm{Spin}_7$, as well as all their Riemannian products. In the non-simply-connected case, examples include flat tori $\mathbb{T}^n$ endowed with the trivial spin structure and their Riemannian products with the above simply connected examples (see also \cite{WangMcK2}).
	\end{example}
	
	The existence of a charged parallel spinor has strong geometric consequences. The first ones are summarized in the following proposition.
	\begin{proposition}\label{PropCPS-1}
	    Let $(M^n,g,E)$ be a charged spin manifold carrying a negative charged parallel spinor $\Phi_-\in\Gamma(S_g)$ (resp. a positive charged parallel spinor $\Phi_+\in\Gamma(S_g)$). Then $\Phi_-$ (resp. $\Phi_+$) has no zero on $M$, the $1$-form associated to the electric field $E$ is given by
		\begin{equation}\label{ElectricFieldSpinor}
			E^\flat=-\frac{1}{n-2}\,d\ln V
			\qquad (\text{resp.}\; E^\flat=\frac{1}{n-2}\,d\ln V)
		\end{equation}
		where $V:=|\Phi_-|^2$ (resp. $V:=|\Phi_+|^2$) and  
		\begin{equation}\label{Ric-CPS}
			\begin{aligned}
				{\rm Ric} 
				&= (n-2)\big(|E|^2g - E^\flat\otimes E^\flat - \nabla E^\flat\big) - \div(E)g,\\
				(\text{resp.}\quad
				{\rm Ric} 
				&= (n-2)\big(|E|^2g - E^\flat\otimes E^\flat + \nabla E^\flat\big) + \div(E)g).
			\end{aligned}
		\end{equation}
	\end{proposition}
	
	\pf
	First, note that for all $X\in\Gamma(TM)$ we have
	\[
	X(V) 
	= 
	2{\rm Re\,}\<\nabla_X\Phi_-,\Phi_-\> 
	= 
	-
	{\rm Re\,}\<c(X)c(E)\Phi_-,\Phi_-\>
	-
	(n-1)g(E,X)|\Phi_-|^2.
	\]
	On the other hand, since 
	\[
	{\rm Re\,}\<c(X)c(E)\Phi,\Phi\>=-g(E,X)|\Phi_-|^2,
	\] 
	we deduce that
	\[
	X(V) 
	= 
	-(n-2)g(E,X)V
	\]
	for all $X\in\Gamma(TM)$. Since $\Phi_-$ is non-trivial, this implies that $\Phi_-$ is nowhere vanishing on $M$, so that $V>0$ and~\eqref{ElectricFieldSpinor} follows. Standard computations show that the curvature of the spinor bundle $S_g$, which by definition is
	\[
	\mathcal{R}_{X,Y}\varphi:=\big([\nabla_X,\nabla_Y]-\nabla_{[X,Y]}\big)\varphi
	\]
	for $\varphi\in\Gamma(S_g)$ and $X$, $Y\in\Gamma(TM)$, satisfies
	\begin{align*}
		\mathcal{R}_{X,Y}\Phi _-
		&= -\frac{1}{2}\big(c(Y)c(\nabla_XE)-c(X)c(\nabla_Y E)\big)\Phi_- \\
		&\quad -\frac{n-1}{2}\big(g(\nabla_XE,Y)-g(\nabla_YE,X)\big)\Phi_- \\
		&\quad -\frac{1}{2}\big(g(E,X)c(Y)-g(E,Y)c(X)\big)c(E)\Phi_- \\
		&\quad -\frac{1}{4}|E|^2[c(X),c(Y)]\Phi_-
	\end{align*}
	for all $X,Y\in\Gamma(TM)$. Now since $E^\flat$ is closed by~\eqref{ElectricFieldSpinor}, we have
	\[
	g(\nabla_XE,Y)=g(\nabla_YE,X)
	\]
	for all $X,Y\in\Gamma(TM)$, so that the above identity simplifies to
	\begin{align*}
		\mathcal{R}_{X,Y}\Phi_-
		&= -\frac{1}{2}\big(c(Y)c(\nabla_XE)-c(X)c(\nabla_Y E)\big)\Phi_- \\
		&\quad -\frac{1}{2}\big(g(E,X)c(Y)-g(E,Y)c(X)\big)c(E)\Phi_- \\
		&\quad -\frac{1}{4}|E|^2[c(X),c(Y)]\Phi_-.
	\end{align*}
	Now let $(e_1,\dots,e_n)$ be a local orthonormal frame of $TM$. The Ricci identity (see \cite{BourguignonHijaziMilhoratMoroianu}) 
	\[
	\sum_{j=1}^n c(e_j)\mathcal{R}_{X,e_j}\Phi_-
	=
	-\frac{1}{2}c\left({\rm Ric}(X)\right)\Phi_-
	\]
	and classical computations, again using that $E^\flat$ is closed, lead to the relation
	\[
	(n-2)c\left(\left(g(E,X)E-|E|^2X\right)+\nabla_XE+\frac{1}{n-2}\div(E)X\right)\Phi_-=-c\left({\rm Ric}(X)\right)\Phi_-
	\]
	and so~\eqref{Ric-CPS} follows since $\Phi_-$ has no zero. The case of positive charged parallel spinors follows by replacing $E$ with $-E$.
	\qed
	
	We now continue the study of charged parallel spinors by analyzing their squared norm.
	\begin{proposition}\label{HessianNormPhi}
		Let $(M^n,g,E)$ be a charged spin manifold carrying a negative charged parallel spinor $\Phi_-$ (resp. a positive charged parallel spinor $\Phi_+$), and set $V:=|\Phi_-|^2$ (resp. $V:=|\Phi_+|^2$). Then $d(VE^\flat)=0$ and 
		\[
		\nabla^2V-(\Delta V)g-V{\rm Ric}
		=
		(n-1)(n-2)V(E^\flat\otimes E^\flat-|E|^2g)
		+(n-1)V\div(E)\,g
		\]
		(resp.
		\[
		\nabla^2V-(\Delta V)g-V{\rm Ric}
		=
		(n-1)(n-2)V(E^\flat\otimes E^\flat-|E|^2g)
		-(n-1)V\div(E)\,g).
		\]
		
	\end{proposition}
	\pf
	It follows from~\eqref{ElectricFieldSpinor} that 
	\[
    VE^\flat=-\frac{1}{n-2}dV
	\]
	and therefore $d(VE^\flat)=0$. On the other hand, it is immediate to check that
	\[
	\nabla^2V
	=
	\nabla dV
	=
	-(n-2)(dV\otimes E^\flat+V\nabla E^\flat)
	\]
	so that 
	\begin{equation}\label{HessianV}
		\nabla^2V 
		=
		(n-2)V\left((n-2)E^\flat\otimes E^\flat-\nabla E^\flat\right).
	\end{equation}
	Taking the trace of this identity proves that
	\begin{equation}\label{LaplacianV}
		\Delta V
		=
		(n-2)V\left((n-2)|E|^2-\div(E)\right). 
	\end{equation}
	Finally putting together~\eqref{Ric-CPS},~\eqref{HessianV} and~\eqref{LaplacianV} yields the second identity. Once again, the case of positive charged parallel spinors follows by replacing $E$ with $-E$.
	\qed 
	
	As a direct consequence of Proposition~\ref{HessianNormPhi}, we obtain that on closed manifolds with $\div(E)=0$, charged parallel spinors reduce to classical parallel spinors.
	\begin{corollary}
		A charged parallel spinor field $\Phi\in\Gamma(S_g)$ on a closed charged spin manifold $(M^n,g,E)$ with $\div(E)=0$ is necessarily parallel. 
	\end{corollary}
	
	\pf 
	It follows from (\ref{LaplacianV}) that
	\[
	\Delta V=(n-2)^2|E|^2V.
	\]
	Integrating over $M$ and using Green's formula gives
	\[
	0=\int_M\Delta V d\mu=(n-2)^2\int_M|E|^2V d\mu.
	\]
	Since $V>0$, we conclude that $E\equiv0$. Hence $\Phi$ is parallel.
	\qed
	
	We now show that charged parallel spinors naturally give rise to electrostatic Einstein--Maxwell systems. The identities obtained in Proposition~\ref{HessianNormPhi} show that $(M^n,g,V,E)$ satisfies a system of equations of electrostatic type. We therefore recall the notion of an electrostatic system (see \cite{CruzLimaDeSousa,FreitasLeandroRibeiroSabo}).
	
	A quadruple $(M^n,g,V,E)$, where $(M^n,g,E)$ is a charged manifold and $V\in C^\infty(M)$, is called an \emph{electrostatic system} if
	\[
	\left\{
	\begin{array}{rll}
		\nabla^2V-(\Delta V)g-V{\rm Ric} &=& (n-1)(n-2)V(E^\flat\otimes E^\flat-|E|^2g),\\
		d(VE^\flat)=0,\,\,\div(E) &=& 0.
	\end{array}
	\right.
	\]
	When these equations are satisfied, the Lorentzian metric
	\[
	\mathbf g=-V^2dt^2+g
	\]
	together with the Maxwell field
	\[
	\mathbf F=\sqrt{\frac{(n-1)(n-2)}{2}}\,dt\wedge E^\flat
	\]
	defines a static solution of the Einstein--Maxwell equations. Comparing with Proposition~\ref{HessianNormPhi}, we immediately obtain the following consequence.
	\begin{corollary}
		Let $(M^n,g,E)$ be a charged spin manifold carrying a charged parallel spinor $\Phi$, and set $V:=|\Phi|^2$. If $\div(E)=0$, then $(M^n,g,V,E)$ is an electrostatic system.
	\end{corollary}
	
	In Section~\ref{MP-Example}, we will see that a particularly important example of electrostatic systems for our purposes is given by the Majumdar--Papapetrou solutions. Their construction from harmonic functions on Euclidean space suggests a close relation between charged parallel spinors and conformal deformations of metrics carrying ordinary parallel spinors. The following result makes this relation precise.
	\begin{proposition}\label{Conformal-PCS}
		A charged spin manifold $(M^n,g,E)$ carries a negative (resp. positive) charged parallel spinor $\Phi_-\in\Gamma(S_g)$ (resp. $\Phi_+\in\Gamma(S_g)$) if and only if there exists a positive smooth function $V$ on $M$ such that 
		\[
		E^\flat=-\frac{1}{n-2}\,d\ln V
		\qquad (\text{resp.}\; E^\flat=\frac{1}{n-2}\,d\ln V)
		\]
		and the conformally related metric $\overline{g}=V^{2/(n-2)}g$ carries a parallel spinor. In this case, we have $V=|\Phi_-|^2$ (resp. $V=|\Phi_+|^2$).
	\end{proposition}
	
	\pf
	Let $\Phi_-\in\Gamma(S_g)$ be a negative charged parallel spinor. Set
	\[
	V:=|\Phi_-|^{2}\ \text{ and }\ \overline{g}=V^{\frac{2}{n-2}}g
	\]
	which is well-defined since $\Phi_-$ is nowhere vanishing. The relation between $E^\flat$ and $V$ follows directly from~\eqref{ElectricFieldSpinor}.
	
	As recalled in Section~\ref{Spinors}, the spinor bundles $S_g$ and $S_{\overline g}$ can be identified via the bundle isomorphism $\beta$ for any metric $\overline{g}\in[g]$. Using~\eqref{ConformalLCC}, we obtain for every $X\in\Gamma(TM)$ and $\varphi\in\Gamma(S_g)$ that
	\[
	\overline{\nabla}_X\beta\varphi
	=
	\beta\Big(\nabla_X\varphi+\frac12 c(X)c(E)\varphi+\frac12 g(E,X)\varphi\Big).
	\]
	Applying this identity to $\Psi_-:=V^{-1/2}\Phi_-$, we obtain
	\[
	\overline{\nabla}_X\beta\Psi_- = \beta\left(X(V^{-1/2})\Phi_-+V^{-1/2}\nabla_X\Phi_-+\frac{1}{2}V^{-1/2}c(X)c(E)\Phi_-+\frac{1}{2}V^{-1/2}g(E,X)\Phi_-\right).
	\]
	Since $\Phi_-$ is a negative charged parallel spinor, it follows from \eqref{ElectricFieldSpinor} that the above expression vanishes, hence $\beta\Psi_-\in\Gamma(S_{\overline g})$ is parallel with respect to $\overline{g}$. The statement for positive charged parallel spinors is obtained by applying the previous argument to the charged manifold $(M^n,g,-E)$.
	
	Conversely, assume that $\overline{g}=V^{2/(n-2)}g$ carries a parallel spinor $\beta\Psi$. Define
	\[
	\Phi:=V^{\frac{1}{2}}\Psi\in\Gamma(S_g)
	\quad
	\text{ and }\quad E^\flat:=-\frac{1}{n-2}d\,\ln V \qquad (\text{resp.}\; E^\flat=\frac{1}{n-2}\,d\ln V).
	\]
	Applying~\eqref{ConformalLCC}, it follows that $\Phi$ is a negative (resp. positive) charged parallel spinor on $(M^n,g,E)$.
	\qed
	
	The conformal characterization also leads to the following observation. If $(M^n,g,E)$ carries a negative charged parallel spinor, then~\eqref{LaplacianV} can be rewritten as
	\begin{equation}\label{EquationU-gbar}
		\Delta_{\bar g}U
		=
		(n-2)\,\div(E)\,U^{\frac{n}{n-2}},
	\end{equation}
	where $\bar g = V^{\frac{2}{n-2}}g$ is the conformal metric appearing in Proposition~\ref{Conformal-PCS} and $U=V^{-1}$. In particular, the Maxwell equation $\div(E)=0$ is equivalent to the harmonicity of $U$ with respect to $\bar g$. This observation suggests a natural way to construct electrostatic systems carrying negative charged parallel spinors. Starting from a Riemannian spin manifold carrying a parallel spinor, any positive harmonic function gives rise, through the above conformal ansatz, to a charged spin manifold admitting a negative charged parallel spinor. A natural source of such examples is provided by Bär's correspondence between real Killing spinors and parallel spinors on metric cones \cite{Bar2}. More precisely, if $(\Sigma^{n-1},h)$ carries a real Killing spinor, then the metric cone
	\[
	(\mathcal{C}(\Sigma),g_C)=((0,\infty)\times \Sigma,\; dr^2+r^2 h)
	\]
	admits a non-trivial parallel spinor, and the above conformal construction can be applied. Let $U$ be a positive harmonic function on $(\mathcal{C}(\Sigma),g_C)$ and define
	\[
	g=U^{\frac{2}{n-2}}g_C,
	\qquad
	E^\flat=\frac{1}{n-2}d\ln U.
	\]
	By Proposition~\ref{Conformal-PCS}, the charged manifold $(M^n,g,E)$ carries a negative charged parallel spinor. A natural class of such harmonic functions is given by radial solutions
	\[
	U(r)=1+\frac{m}{r^{n-2}},
	\]
	where $m>0$, which are harmonic on the cone. Near the singular point $\{r=0\}$, the metric $g$ is asymptotically cylindrical. More precisely, setting \(t=-m^{\frac1{n-2}}\ln r\), one has
	\[
	g\sim dt^2+m^{\frac{2}{n-2}}h\quad \text{as}\quad t\to+\infty
	\]
	so that the end is asymptotically cylindrical with limiting cross-section \((\Sigma^{n-1},m^{\frac{2}{n-2}}h)\). At infinity, we have \(V\to 1\), so that the metric $g$ is asymptotic to $g_C$. In particular, the metric is asymptotically flat if and only if $(\Sigma^{n-1},h)$ is isometric, up to constant scaling, to the round sphere, in which case one recovers the extremal Reissner--Nordstr\"om manifold (see Section~\ref{MP-Example}). 
	
	The previous results were obtained under the existence of a single charged parallel spinor. It is natural to ask what happens when a charged manifold carries both a positive and a negative charged parallel spinor. The following theorem provides a complete description of this situation.
	
	\begin{theorem}
		Let $(M^n,g,E)$ be a complete charged spin manifold with $E\not\equiv 0$. Then the following are equivalent:
		\begin{enumerate}
			\item $(M^n,g,E)$ carries both a positive and a negative charged parallel spinor;
			\item $(M^n,g)$ is isometric to a Riemannian product
			\[
			(\mathbb{R}\times \Sigma,dt^2+h),
			\]
			where $(\Sigma^{n-1},h)$ is a complete Riemannian spin manifold admitting a real Killing spinor with Killing constant $\kappa$ and $E=2\kappa\partial_t$.
		\end{enumerate}
	\end{theorem}
	
	We first recall that a spinor field $\psi\in\Gamma(S_h)$ on $(\Sigma^{n-1},h)$ is called a \emph{real Killing spinor with Killing constant $\kappa\in\mathbb{R}\setminus\{0\}$} if it satisfies the overdetermined equation
	\[
	\nabla^h_X\psi=\kappa c_h(X)\psi	
	\]
	for all $X\in\Gamma(T\Sigma)$. It is well known that the existence of such a spinor field on a complete Riemannian spin manifold forces the metric to be Einstein with positive scalar curvature. In particular, such a manifold is compact with finite fundamental group. In the simply connected case, complete manifolds carrying real Killing spinors were classified by B\"ar~\cite{Bar2}. Besides the round spheres, this includes Einstein--Sasaki manifolds, $3$-Sasakian manifolds, and, in dimensions $6$ and $7$, nearly K\"ahler manifolds and nearly parallel $G_2$ manifolds.
	
	\vspace{0.5cm}
	
	\pf
	Assume first that $(M^n,g,E)$ carries both a positive and a negative charged parallel spinor. Then both identities in~\eqref{Ric-CPS} hold simultaneously. Subtracting them, we obtain
	\[
	\nabla_XE=-\frac{1}{n-2}\div(E) X
	\]
	for all $X\in\Gamma(TM)$. Taking the trace gives
	\[
	\div(E)=-\frac{n}{n-2}\operatorname{div}(E),
	\]
	hence $\div(E)=0$, and therefore $\nabla E=0$. On the other hand, it follows from~\eqref{ElectricFieldSpinor} that
	\[
	E=\frac{1}{n-2}\nabla\ln V,
	\qquad
	V:=|\Phi_+|^2
	\]
	where $\Phi_+\in\Gamma(S_g)$ denotes a positive charged parallel spinor. Since $E$ is parallel, the function $\ln V$ is affine. Therefore, by the splitting theorem of Innami~\cite{Innami}, $(M^n,g)$ is isometric to a Riemannian product
	\[
	\left(\mathbb{R}\times \Sigma,dt^2+h\right),
	\]
	where $(\Sigma^{n-1},h)$ is a $(n-1)$-dimensional complete Riemannian manifold. Moreover, after possibly reversing the $t$-coordinate, we may assume that
	\[
	E=|E|\,\frac{\partial}{\partial t}
	\qquad\text{and}\qquad
	V(t)=Ce^{(n-2)|E|t}
	\]
	for some constant $C>0$.
	
	Since $M$ is spin, the hypersurface $\Sigma$ inherits a spin structure. Moreover, a spinor field on $\mathbb{R}\times \Sigma$ may be viewed as a smooth one-parameter family, parametrized by $t\in\mathbb{R}$, of sections of the restricted spinor bundle $\RSB_h$. In particular, it follows from~\eqref{PositiveCPS} that $\Phi_{+}$ satisfies 
	\[
	\nabla_{\frac{\partial}{\partial t}}\Phi_+=\frac{\partial\Phi_+}{\partial t}=\frac{n-2}{2}|E|\Phi_+.
	\]
	Hence, for each fixed $x\in \Sigma$, the map $t\mapsto \Phi_+(t,x)$ satisfies a linear ordinary differential equation in the finite-dimensional vector space $\RSB_{h,x}$. By uniqueness of solutions, there exists $\psi_+\in\Gamma(\RSB_h)$ such that
	\[
	\Phi_+(t,x)=e^{\frac{n-2}{2}|E|t}\psi_+(x)
	\]
	for all $(t,x)\in\mathbb{R}\times\Sigma$. Since $\Sigma$ is a totally geodesic hypersurface of $(M^n,g)$, it follows from~\eqref{CliffordLCI-Extrinsic} that the extrinsic spin Levi-Civita connection $\nb^h$ coincides with the ambient spin connection $\nabla$ along $\Sigma$ so that~\eqref{PositiveCPS} yields
	\[
	\nb^h_X\psi_+
	=
	\frac{|E|}{2}\,c(X)c(\nu)\psi_+
	\]
	for all $X\in\Gamma(T\Sigma)$ and where $\nu:=\partial/\partial t$ is a unit normal to $\Sigma$ in $M$. Thus $\psi_+$ is an extrinsic Killing spinor on $(\Sigma^{n-1},h)$ in the sense of~\cite{HijaziMontiel1}. Using the identifications~\eqref{SpinorBundleRestriction} and~\eqref{CliffordLCI-Identifications}, it follows that $\psi_+$ induces on $(\Sigma^{n-1},h)$ a real Killing spinor with Killing constant $|E|/2$ or $-|E|/2$.
	
	Conversely, assume that $(M^n,g)$ is isometric to a Riemannian product
	\[
	(\mathbb{R}\times \Sigma, dt^2+h),
	\]
	where $(\Sigma^{n-1},h)$ is a complete Riemannian spin manifold carrying a real Killing spinor $\psi$ with Killing constant $\kappa\neq 0$, and that $E=2\kappa\,\partial_t$. We show that $(M^n,g,E)$ carries both a positive and a negative charged parallel spinor. We first observe that $M$ is naturally endowed with a spin structure and, using the identifications~\eqref{SpinorBundleRestriction} and~\eqref{CliffordLCI-Identifications}, the real Killing spinor $\psi\in\Gamma(S_h)$ gives rise to an extrinsic Killing spinor $\psi_+\in\Gamma(\RSB_h)$ with Killing constant $\kappa$. More precisely,
	\[
	\psi_+ :=
	\begin{cases}
		\psi & \text{if } n \text{ is odd},\\
		(\psi,0) & \text{if } n \text{ is even}.
	\end{cases}
	\]
	Then, it follows from~\eqref{CliffordLCI-Extrinsic} that
	\[
	\psi_-:=c(\nu)\psi_+
	\]
	is an extrinsic Killing spinor with Killing constant $-\kappa$. Without loss of generality, we may assume that $\kappa>0$. Now define
	\[
	\Phi_+:=e^{(n-2)\kappa t}\psi_+\quad\text{and}\quad\Phi_-:=e^{-(n-2)\kappa t}\psi_-.
	\]
	A straightforward computation, similar to that of the direct implication, shows that $\Phi_+$ and $\Phi_-$ define respectively a positive and a negative charged parallel spinor on $(M^n,g,E)$.
	\qed
	
	%%%%%%%%%%%%%%%%%%%%%%%%%%%%%%%%%%%%%%%%%%%%%%%%%%%%%%%%%%%%%%%%%%%%%%%%%%%%
	
	\section{The Majumdar-Papapetrou manifolds}\label{MP-Example}
	
	%%%%%%%%%%%%%%%%%%%%%%%%%%%%%%%%%%%%%%%%%%%%%%%%%%%%%%%%%%%%%%%%%%%%%%%%%%%%
	
	As observed at the end of the previous section, the construction of charged parallel spinors from parallel spinors on cones naturally gives rise to the extremal Reissner--Nordstr\"om manifold. It is therefore natural to place this example in the broader context of the Majumdar--Papapetrou family, in which the extremal Reissner--Nordstr\"om manifold appears as a particular case.
	
	Let $x_1,\dots,x_\ell\in\mathbb{R}^n$, $m_1,\dots,m_\ell\in\mathbb{R}_+^*$ and $\ell\ge1$. Set
	\[
	M^n=\mathbb{R}^n\setminus\{x_1,\dots,x_\ell\}\quad\text{and}\quad m:=\sum_{j=1}^\ell m_j>0,
	\]
	and consider the metric
	\[
	g_{m}=U_{m}^{\frac{2}{n-2}}\delta,
	\]
	where
	\[
	U_{m}(x)=1+\sum_{j=1}^\ell \frac{m_j}{|x-x_j|^{n-2}}.
	\]
	The electric field is defined by
	\[
	E_{m}^\flat=-\frac{1}{n-2}d\ln V_{m},
	\qquad
	V_m=U_m^{-1}.
	\]
	Then $(M^n,g_m,E_m)$ defines an electrostatic system in the sense of Section~\ref{Section-CPS}. Moreover, $(M^n,g_m)$ has one asymptotically flat end and $\ell$ asymptotically cylindrical ends (as defined in Section~\ref{MainDefinitions}), corresponding to the region at infinity and to the neighbourhood of each point $x_j$. The limiting cross-section of the cylindrical end associated with $x_j$ is
	\[
	(\mathbb S^{n-1},\,m_j^{\frac{2}{n-2}}g_{\mathbb S^{n-1}}).
	\]
	We will refer to the charged manifold $(M^n,g_m,E_m)$ as a \emph{Majumdar--Papapetrou manifold of mass $m$}. With the above convention for the electric field, the total charge of such an asymptotically flat manifold is $Q=-m <0$. Replacing $E$ by $-E$ produces again an electrostatic system with the same metric, but with opposite charge. In both cases, the ADM mass satisfies $m=|Q|$. 
	
	From the perspective of spin geometry, Proposition~\ref{Conformal-PCS} applies to a Majumdar--Pa\-pa\-pe\-trou manifold. Indeed, with the above choice of the electric field, the conformal metric $V_m^{2/(n-2)}g_{m}$ is the Euclidean metric $\delta$. In particular, it carries a maximal number of parallel spinors, which gives rise to a maximal family of charged parallel spinors on a Majumdar--Papapetrou manifold. More precisely, if $\{\phi_j\,|\, j=1,\dots,N\}$ is an orthonormal basis of parallel spinors of $(\mathbb{R}^n,\delta)$, then
	\[
	\Phi_j := V^{1/2}_m\beta\phi_j,\,\quad j=1,\dots,N
	\]
	defines a basis of charged parallel spinors; moreover, $|\Phi_j|^2 = V_m$ for all $j=1,\dots,N$. Here $N=2^{[n/2]}$ is the rank of the (complex) spinor bundle, and $\beta$ denotes the identification between the spinor bundles over $M$ with respect to \(\delta\) and \(g_m\) (see Section~\ref{Spinors}). With our choice of sign convention, these spinors are necessarily \emph{negative}. This is consistent with the fact that the total charge of the asymptotically flat end is negative. Replacing $E_m$ by $-E_m$ produces a positively charged manifold, and the corresponding charged parallel spinors become positive.
	
	At infinity, we have $V_m\to 1$, so that $\Phi_j$ tends to a constant spinor. Near each pole $x_j$, writing $r_j=|x-x_j|$ and \(t=-m^{1/(n-2)}_j\ln r_j\), one has
	\[
	g_m\sim dt^2+ m_j^{\frac{2}{n-2}}g_{\mathbb S^{n-1}}
	\]
	and
	\begin{equation}\label{CPS-CylindricalBehaviour}
		|\Phi_j|\sim \frac{1}{\sqrt{m_j}}e^{-\frac{n-2}{2}m_j^{-\frac{1}{n-2}}t}\quad \text{as } t\to+\infty.
	\end{equation}
	In particular, the associated charged parallel spinors decay exponentially along each cylindrical end. We summarize the above discussion in the following proposition.
	\begin{proposition}
		A Majumdar-Papapetrou manifold admits the maximal number of charged parallel spinors. These spinors are asymptotic to constant spinors at infinity and satisfy~\eqref{CPS-CylindricalBehaviour} along each cylindrical end.
	\end{proposition}
	
	For the boundary rigidity statement, we will need exterior regions in the model geometries with compact inner boundary. We therefore truncate
	the cylindrical ends of a Majumdar--Papapetrou manifold as follows. For $r_0>0$, we define
	\[
	\Omega_{m,r_0}:= \mathbb{R}^n \setminus \bigcup_{j=1}^\ell B_{r_0}(x_j),
	\]
	where $B_{r_0}(x_j)$ denotes the Euclidean ball of radius $r_0$ centered at $x_j$. We view $\Omega_{m,r_0}$ as a domain in a Majumdar--Papapetrou manifold of mass $m$. Its boundary is given by
	\[
	\partial\Omega_{m,r_0} = \bigsqcup_{j=1}^\ell \partial\Omega_{m,r_0}^{(j)},
	\qquad
	\partial\Omega_{m,r_0}^{(j)}:=\{x\in\mathbb{R}^n \mid |x-x_j|=r_0\}.
	\]
	We denote by $\nu$ the unit normal vector field along $\partial\Omega_{m,r_0}$ pointing toward infinity. Using the fact that the metric $g_m$ is conformally flat, it is straightforward to compute that the Weingarten map of a connected component of $\partial \Omega_{m,r_0}$ satisfies
	\begin{equation}\label{MPSphere-WP}
		A^{(j)}_mX=\big(r_0^{-1}V_m^{\frac{1}{n-2}}+g_m(E_m,\nu)\big)X
	\end{equation}
	for all $X\in\Gamma(T\partial\Omega_{m,r_0}^{(j)})$ so that it is a totally umbilical hypersurface with mean curvature 
	\[
	H^{(j)}_m=(n-1)\big(r_0^{-1}V_m^{\frac{1}{n-2}}+g_m(E_m,\nu)\big).
	\]
    The charged parallel spinors constructed above satisfy a remarkably simple boundary equation on these hypersurfaces.
    \begin{proposition}\label{CPS-RestrictionMP}
    	Let $\Phi$ be a negative charged parallel spinor on $(M^n,g_m,E_m)$. For each $j=1,\dots,\ell$, let $\D^{(j)}$ denote the extrinsic Dirac operator of
    	$\partial\Omega_{m,r_0}^{(j)}$ endowed with the induced metric. Then
    	\[
    	\D^{(j)}\Phi
    	=
    	-\frac{n-1}{2}
    	r_0^{-1}V_m^{\frac1{n-2}}\Phi
    	\]
    	along $\partial\Omega_{m,r_0}^{(j)}$.
    \end{proposition}

	\pf 
	Let \(\nb^{(j)}\) denote the extrinsic spin connection of \(\partial\Omega_{m,r_0}^{(j)}\). Writing 
	\[
	E_m=E_m^\top+g_m(E_m,\nu)\nu
	\]
	where $E_m^\top$ denotes the tangential component of \(E_m\) along \(\partial\Omega_{m,r_0}^{(j)}\), it follows from the definition of a negative charged parallel spinor,~\eqref{CliffordLCI-Extrinsic} and~\eqref{MPSphere-WP} that
	\[
	\nb^{(j)}_X\Phi=\frac{1}{2}r_0^{-1}V_m^{\frac{1}{n-2}}c(X)c(\nu)\Phi-\frac{1}{2}c(X)c(E^\top_m)\Phi-\frac{n-1}{2}g_m(E^\top_m,X)\Phi
	\]
	for all $X\in\Gamma(T\partial\Omega_{m,r_0}^{(j)})$. Then from the definition~\eqref{ExtrinsicDirac} of the extrinsic Dirac operator, a direct computation leads to the announced formula. 
	\qed
	
	In the special case $\ell=1$, the Majumdar--Papapetrou manifold of mass $m$ is the so-called \emph{extremal Reissner--Nordstr\"om manifold of mass $m$}. Here, without loss of generality, one can assume that $x_1=0$ and, setting $r=|x|$, we have
	\[
	U_{m}(x)=1+\frac{m}{r^{n-2}}.
	\]
	In particular, both \(U_m\) and \(V_m\) are radial functions, and the electric field is purely radial. This additional symmetry is absent in the general Majumdar--Papapetrou case \((\ell\geq 2)\). For \(r_0>0\), the domain \(\Omega_{m,r_0}\) is then an asymptotically flat manifold whose boundary is a round sphere with metric denoted by \(h_m\), and is totally umbilical with constant mean curvature
	\[
	H_m=(n-1)\left((m+r_0^{n-2})^{-\frac{1}{n-2}}+g_m(E_m,\nu)\right).
	\]
	Since \(V_m\) is constant along \(\partial\Omega_{m,r_0}\), the boundary equation of Proposition~\ref{CPS-RestrictionMP} reduces to the extrinsic Killing spinor equation. 
	\begin{corollary}
		Let $\Phi$ be a negative charged parallel spinor on the extremal Reissner-Nordstr\"om manifold of mass $m$. Its restriction to $\partial\Omega_{m,r_0}$ is an extrinsic Killing spinor with a positive Killing constant. More precisely, it satisfies
		\[
		\nb^{h_m}_X\Phi=\frac{1}{2}(m+r_0^{n-2})^{-\frac{1}{n-2}}c(X)c(\nu)\Phi
		\]
		for all $X\in \Gamma(T\partial\Omega_{m,r_0})$ and in particular
		\[
		\D_{h_m}\Phi=-\frac{n-1}{2}(m+r_0^{n-2})^{-\frac{1}{n-2}}\Phi.
		\]
	\end{corollary}
	In particular, negative charged parallel spinors on the extremal Reissner--Nordstr\"om manifold induce boundary data satisfying the Atiyah--Patodi--Singer condition. This observation reveals the geometric structure underlying the equality case in Theorems~\ref{ChargedPMT-Connected boundary} and~\ref{ChargedPMT-Connected-3dim}, and explains why the extremal Reissner--Nordstr\"om manifold naturally emerges as the rigidity model in the results proved below.
	
	%%%%%%%%%%%%%%%%%%%%%%%%%%%%%%%%%%%%%%%%%%%%%%%%%%%%%%%%%%%%%%%%%%%%%%%%%%%
	
	\section{Manifolds with compact boundary}\label{PMT-CB}
	
	%%%%%%%%%%%%%%%%%%%%%%%%%%%%%%%%%%%%%%%%%%%%%%%%%%%%%%%%%%%%%%%%%%%%%%%%%%%
	
	We now turn to the first geometric setting in which a non-trivial electric charge may occur, namely asymptotically flat manifolds with compact inner boundary.
	
	Throughout this section, we consider complete charged spin manifolds $(M^n,g,E)$ with compact boundary and a single asymptotically flat end. As explained at the end of the proof of Theorem~\ref{ChargedPMT-Connected boundary}, the arguments extend without difficulty to the case of finitely many asymptotically flat ends.
	
	%%%%%%%%%%%%%%%%%%%%%%%%%%%%%%%%%%%%%%%%%%%%%%%%%%%%%%%%%%%%%%%%%%%%%%%%%%%
	
	\subsection{Solving the modified Dirac equation}
	
	%%%%%%%%%%%%%%%%%%%%%%%%%%%%%%%%%%%%%%%%%%%%%%%%%%%%%%%%%%%%%%%%%%%%%%%%%%%
	
	We establish the analytic framework required for the proof of Theorem~\ref{ChargedPMT-Connected boundary}. We first recall the Atiyah--Patodi--Singer (APS) boundary condition. Let \(h:=g_{|\partial M}\) denote the induced metric on \(\partial M\). The extrinsic Dirac operator $\D_h$ is a self-adjoint elliptic operator whose spectrum consists of a discrete unbounded sequence of real eigenvalues $(\lambda_k)_{k\in\mathbb{Z}}$. Moreover, the relation
	\begin{equation}\label{ExtrinsicDiracSpectrum}
		\D_h\big(c(\nu)\varphi\big)=-c(\nu)\D_h\varphi
	\end{equation}
	implies that the spectrum of $\D_h$ is symmetric with respect to zero. Accordingly, we may index the eigenvalues so that
	\[
	\lambda_k>0 \text{ for } k>0, \qquad \lambda_k=-\lambda_{-k} \text{ for } k<0, \qquad \lambda_0=0.
	\]
	In the sequel, we assume that $\D_h$ is invertible. We define the APS projection
	\[
	P_+ : L^2(\RSB_h)\longrightarrow L^2(\RSB_h)
	\]
	as the $L^2$-orthogonal projection onto the subspace spanned by the eigenspinors associated with positive eigenvalues of $\D_h$. It is well known (see, for instance,~\cite{BartnikChrusciel}) that $P_+$ defines a global elliptic boundary condition for Dirac-type operators such as $D^\mp$.
	
	Using Remark~\ref{IPP-DM} and~\eqref{ExtrinsicDiracSpectrum}, one checks that the formal adjoint of $D^\mp$ under the boundary condition $P_+$ is $D^\pm$ equipped with the same boundary condition.
	
	Theorem~8.9 of Bartnik and Chru\'sciel~\cite{BartnikChrusciel} requires a weighted Poincar\'e inequality together with a Schr\"odinger--Lichnerowicz estimate. We verify below that both assumptions hold in our setting.

	\medskip
	
	$\bullet$ {\it A weighted Poincar\'e inequality}, that is, there exists $w\in L^1_{loc}(M)$ with ${\rm ess\,inf}_\Omega w>0$ for every relatively compact domain $\Omega \subset M$ such that for all $\varphi\in C^1_c(M,S_g)$, the space of $C^1$ compactly supported spinor fields on $M$, we have 
	\begin{equation}\label{wPi}
		\int_M|\varphi|^2w\,d\mu\leq\int_M |\nabla^\mp\varphi|^2d\mu.
	\end{equation}
	A direct computation shows that $\Gamma^\mp$, the symmetric part of the connections $\nabla^\mp$ is given by
	\[
	\Gamma^\mp_X=\mp\frac{n-2}{2}g(E,X)\in\Gamma\big({\rm End}(S_g)\big).
	\]
	Then our decay assumptions~\eqref{AFDecay} imply that the hypotheses of \cite[Theorem 9.5]{BartnikChrusciel} are satisfied, and hence~\eqref{wPi} follows.
	
	\medskip
	
	$\bullet$ {\it A Schr\"odinger--Lichnerowicz estimate for the pair $(D^\mp,P_+)$}. We first note that the boundary term in the integral formula of Corollary \ref{IntegralVersion} takes the form
	\begin{equation}\label{RestrictedDiracEquation}
		L^{\mp}_\nu\varphi=-\D_h\varphi-\frac{1}{2}\mathcal{H}^\mp\varphi
	\end{equation}
	where 
	\[
	\mathcal{H}^\mp:=H\mp(n-1)g(E,\nu).
	\]
	It follows that
	\begin{equation}\label{IntegralVersionBoundary-1}
		\oint_{\partial M}\<\D_h\varphi+\frac{1}{2}\mathcal{H}^\mp\varphi,\varphi\>d\sigma + \int_M|D^\mp\varphi|^2 d\mu= \int_{M}\big(|\nabla^\mp\varphi|^2+\mathcal{R}^\mp|\varphi|^2\big)d\mu
	\end{equation}
	for all $\varphi\in C^1_c(M,S_g)$. Since $\D_h$ is assumed to be invertible, the trace of $\varphi$ along $\partial M$ admits the spectral decomposition
	\[
	\varphi
	=
	\sum_{k\in\mathbb Z^\ast}A_k\varphi_k
	=
	\sum_{k<0}A_k\varphi_k
	+
	P_+\varphi,
	\]
	where $(A_k)_{k<0}\subset\mathbb C$. Here $(\varphi_k)_{k\in\mathbb{Z}^\ast}$ is a complete $L^2$-orthonormal basis of $\RSB_h$ consisting of smooth eigenspinors of the extrinsic Dirac operator $\D_h$ associated with eigenvalues $\lambda_k$. Using this decomposition, one readily checks that
	\begin{align*}
		\oint_{\partial M}\<\D_h\varphi+\frac{1}{2}\mathcal{H}^\mp\varphi,\varphi\>\,d\sigma 
		& \leq \left(1+\frac{1}{2}\sup_{\partial M}\mathcal{H}^\mp\right)\oint_{\partial M}|P_+\varphi|^2\,d\sigma +\left(\frac{1}{2}\sup_{\partial M}\mathcal{H}^\mp - \lambda_1\right)\sum_{k<0}|A_k|^2
	\end{align*}
	for all $\varphi\in C^1_c(M,S_g)$. In particular, if we assume that 
	\[
	\lambda_1 \ge \frac{1}{2}\sup_{\partial M}\left(H+(n-1)|g(E,\nu)|\right),
	\] 
	we obtain
	\[
	\oint_{\partial M}\<\D_h\varphi+\frac{1}{2}\mathcal{H}^\mp\varphi,\varphi\>\,d\sigma \leq  \oint_{\partial M}|\mathcal{J}^\mp P_+\varphi|^2\,d\sigma
	\]
	where $\mathcal{J}^\mp:=\big(|\D_h|+\frac{1}{2}\sup_{\partial M}\mathcal{H}^\mp\big)^{1/2}$. Combining this estimate with the integral identity~\eqref{IntegralVersionBoundary-1} and using that $\mathcal{R}^\mp\ge 0$, we obtain
	\[
	\int_M|\nabla^\mp\varphi|^2\,d\mu\leq \int_M|D^\mp\varphi|^2\,d\mu+\oint_{\partial M}|\mathcal{J}^\mp P_+\varphi|^2\,d\sigma
	\]
	for all $\varphi\in C^1_c(M,S_g)$. This proves that the pair $(D^\mp,P_+)$ satisfies a Schr\"odinger--Lichnerowicz estimate in the sense of~\cite[Definition~8.4]{BartnikChrusciel}.
	
	\medskip
	
	Now consider the space
	\[
	C^1_c(M,S_g;P_+)
	:=
	\{\varphi\in C^1_c(M,S_g)\,;\,P_+\varphi=0\},
	\]
	equipped with the norm
	\[
	\|\varphi\|_{\mp}^2
	:=
	\int_M
	\Big(
	|\nabla^\mp\varphi|^2
	+
	\mathcal R^\mp |\varphi|^2
	\Big)\,d\mu .
	\]
	Since \(\mathcal R^\mp\ge0\) by the dominant energy condition~\eqref{DEC}, this indeed defines a norm. We denote by \(\mathbb H^\mp\) the corresponding Hilbert space completion. Applying Theorem~8.9 of Bartnik and Chru\'sciel~\cite{BartnikChrusciel} now yields the following proposition. 
	\begin{proposition}\label{ExistenceBVP}
		Let $n \geq 3$ and let $(M^n,g,E)$ be a complete charged spin manifold with at least one asymptotically flat end and compact boundary $\partial M$. Assume that
		the dominant energy condition~\eqref{DEC} holds, that the extrinsic Dirac operator has trivial kernel and its first positive eigenvalue satisfies
		\[
		\lambda_1 \ge \frac{1}{2}\sup_{\partial M}\left(H+(n-1)|g(E,\nu)|\right).
		\] 
		Then $D^\mp : \mathbb{H}^\mp \to L^2(S_g)$ is an isomorphism.
	\end{proposition}
	
	%%%%%%%%%%%%%%%%%%%%%%%%%%%%%%%%%%%%%%%%%%%%%%%%%%%%%%%%%%%%%%%%%%%%%%%%%%%%
	
	\subsection{Proof of Theorem~\ref{ChargedPMT-Connected boundary}}\label{ProofTHMCB}
	
	%%%%%%%%%%%%%%%%%%%%%%%%%%%%%%%%%%%%%%%%%%%%%%%%%%%%%%%%%%%%%%%%%%%%%%%%%%%%
	
	Without loss of generality, we assume that $Q\leq 0$, since otherwise one can replace $E$ by $-E$. Let $\Phi_0\in\Gamma(S_g)$ be an asymptotically constant spinor supported in the asymptotically flat end. Under the decay assumptions~\eqref{AFDecay}, we have $D^-\Phi_0\in L^2(S_g)$. The spectral condition required to apply Proposition~\ref{ExistenceBVP} follows from the Bär--Hijazi inequality~\cite{Bar3,Hijazi2,Hijazi1},
	\[
	\lambda_1 \ge \frac{1}{2}\,{\rm Vol}(\partial M,h)^{-\frac{1}{n-1}}
	\sqrt{\frac{n-1}{n-2}\,\mathcal Y(\partial M,[h])},
	\]
	which, combined with~\eqref{YamabeMC}, yields the desired lower bound on $\lambda_1$. Recall that the Yamabe constant of the conformal class \([h]\) is defined by
	\begin{align*}
		\mathcal{Y}(\partial M,[h]):= \inf_{\overline{h}\in[h]} \left\{
		\frac{\oint_{\partial M}\overline{\Ss}\,
			d\sigma_{\overline{h}}}
		{{\rm Vol}(\partial M,\overline{h})^{\frac{n-3}{n-1}}}
		\right\}
	\end{align*}
	where $\overline{\Ss}$ and ${\rm Vol}(\partial M,\overline{h})$ denote respectively the scalar curvature and the volume of $\partial M$ computed with respect to a metric $\overline{h}$ in the conformal class $[h]$. Proposition~\ref{ExistenceBVP} therefore provides a unique $\phi_0\in\mathbb{H}^-$ such that $D^-(\Phi_0+\phi_0)=0$. Setting $\Phi:=\Phi_0+\phi_0$, we obtain a spinor field satisfying
	\[
	D^-\Phi=0 \quad \text{on } M, \qquad P_+\Phi=0 \quad \text{along } \partial M,
	\]
	with $\Phi-\Phi_0\in\mathbb{H}^-$. Applying Corollary~\ref{IntegralVersion} to the region $M_r$ bounded by $\partial M$ and a large coordinate sphere $S_r$, we obtain
	\[
	\oint_{S_r}\langle L^-_{\nu_r}\Phi,\Phi\rangle\,d\sigma_r 
	\;\geq\;
	\int_{M_r}\Big(|\nabla^-\Phi|^2+\mathcal{R}^-|\Phi|^2\Big)\,d\mu
	\;\geq\; 0,
	\]
	because of the dominant energy condition~\eqref{DEC} and the contribution from the inner boundary being nonnegative by~\eqref{YamabeMC}. Since $\Phi-\Phi_0\in\mathbb{H}^-$, it follows that
	\[
	\lim_{r\to\infty}\oint_{S_r}\langle L^-_{\nu_r}\Phi,\Phi\rangle\,d\sigma_r
	=
	\lim_{r\to\infty}\oint_{S_r}\langle L^-_{\nu_r}\Phi_0,\Phi_0\rangle\,d\sigma_r.
	\]
	Moreover,
	\[
	L^-_{\nu_r}\Phi_0
	=
	L_{\nu_r}\Phi_0
	+\frac{n-1}{2}g(E,\nu_r)\,\Phi_0,
	\]
	where $L_{\nu_r}\Phi_0=\nabla_{\nu_r}\Phi_0+c(\nu_r)D\Phi_0$. A standard computation shows that
	\[
	\lim_{r\to\infty}\oint_{S_r}\langle L_{\nu_r}\Phi_0,\Phi_0\rangle\,d\sigma_r
	=
	\frac{n-1}{2}\omega_{n-1}m\,|\Phi_0|^2.
	\]
	On the other hand, since $\Phi_0$ is asymptotically constant,
	\[
	\lim_{r\to\infty}\oint_{S_r}g(E,\nu_r)|\Phi_0|^2\,d\sigma_r
	=
	\omega_{n-1}Q\,|\Phi_0|^2.
	\]
	Substituting these identities into the previous inequality yields $m\ge |Q|$.
	
	We now turn to the rigidity statement and assume that equality holds, that is, $m=|Q|$. Assume moreover that $\div(E)=0$. Returning to the integral identity, equality implies that all inequalities are saturated. In particular, $\Phi$ is a negative charged parallel spinor, i.e. $\nabla^-\Phi=0$. Setting \(V:=|\Phi|^2\), it follows from Proposition~\ref{PropCPS-1} and from $\div(E)=0$ that
	\[
	\Delta_g(\ln V)=0.
	\]
	Moreover, since $\Phi$ is asymptotic to a nonzero constant spinor,
	\[
	\ln V\to 0 \quad \text{as}\quad r\rightarrow\infty.
	\]
	Hence it follows from standard asymptotic expansions for harmonic functions (see~\cite{Lee} for example) that
	\[
	\ln V=\frac{A}{r^{n-2}}+o_2(r^{2-n})
	\]
	for some constant $A\in\mathbb R$. Using this expansion in the expression of $E^\flat$, we obtain
	\[
	E^i = A\frac{x_i}{r^n}+o_1(r^{1-n}),
	\]
	where the components are computed in the asymptotic Euclidean coordinates. This implies that $A=Q$ and in particular,
	\begin{equation}\label{eq:lnV-expansion}
		\ln V=\frac{Q}{r^{n-2}}+o_2(r^{2-n}).
	\end{equation}
	
	Moreover, equality in the boundary estimate forces equality in the Bär--Hijazi inequality. This implies that the restriction of $\Phi$ to $\partial M$ is an extrinsic Killing spinor with Killing number $\lambda_1/(n-1)$. Then it follows from~\eqref{RestrictedDiracEquation} that
	\[
	-\lambda_1\Phi=\D_h\Phi=-\frac{\mathcal{H}^-}{2}\Phi
	\]
	and so $\mathcal{H}^-$ is constant along $\partial M$. Therefore, it holds that
	\[
	\nb^h_X\Phi
	=
	\frac{\mathcal H^-}{2(n-1)}c(X)c(\nu)\Phi
	\]
	for all $X\in\Gamma(T\partial M)$. On the other hand, using the spinorial Gauss formula~\eqref{CliffordLCI-Extrinsic} in the equation $\nabla^-\Phi=0$ yields
	\[
	\nb^h_X\Phi
	=
	-\frac{1}{2}c(X)c(E)\Phi
	-\frac{n-1}{2}g(E,X)\Phi
	+\frac{1}{2}c(AX)c(\nu)\Phi.
	\]
	Comparing both expressions and decomposing $E=E^\top+g(E,\nu)\nu$, we obtain
	\[
	\left(c(X)c(E^\top)+(n-1)g(E^\top,X)\right)\Phi
	=
	c\Big(\frac{1}{n-1}HX-AX\Big)c(\nu)\Phi.
	\]
	Taking the real part of the scalar product with $\Phi$ and using the skew-symmetry of Clifford multiplication, we deduce that $E^\top=0$. It follows that $\ln V$ is constant along $\partial M$, hence $V$ is a positive constant on $\partial M$ since $\partial M$ is connected. Denote by $V_0$ this constant. Substituting back, we obtain
	\[
	\frac{1}{n-1}HX-AX=0,
	\]
	so that $\partial M$ is totally umbilic.
	
	We may therefore apply the conformal characterization established in Section~\ref{Section-CPS}. In particular, Proposition~\ref{Conformal-PCS} shows that the metric
	\[
	\bar g:=V^{\frac{2}{n-2}}g
	\]
	is Ricci-flat. Moreover, the asymptotic expansion~\eqref{eq:lnV-expansion} implies that $\bar g$ is asymptotically flat. On the other hand, the conformal transformation law for the mean curvature gives
	\[
	\bar H
	=
	V_0^{-\frac{1}{n-2}}
	\left(
	H+\frac{n-1}{n-2}\frac{\partial_\nu V}{V_0}
	\right).
	\]
	It follows from~\eqref{ElectricFieldSpinor} that
	\[
	\frac{\partial_\nu V}{V_0}=-(n-2)g(E,\nu),
	\]
	and then
	\[
	\bar H
	=
	V_0^{-\frac{1}{n-2}}\bigl(H-(n-1)g(E,\nu)\bigr)
	=
	V_0^{-\frac{1}{n-2}}\mathcal H^-.
	\]
	Since $\mathcal H^-$ is constant on $\partial M$, the mean curvature $\bar H$ is constant. Furthermore, the boundary condition~\eqref{YamabeMC} translates into
	\[
	\bar H
	\le
	{\rm Vol}(\partial M,\bar{h})^{-\frac{1}{n-1}}
	\sqrt{\frac{n-1}{n-2}\,\mathcal Y(\partial M,[\bar h])}
	\]
	because $V$ is constant on $\partial M$ so that the induced metrics $h$ and $\bar{h}:=\bar g_{|\partial M}$ differ only by a constant scaling factor. Thus $(M^n,\bar g)$ satisfies the hypotheses of Herzlich’s positive mass theorem for manifolds with boundary~\cite{Herzlich1,Herzlich2}. 
	In fact, using the asymptotic expansion~\eqref{eq:lnV-expansion} together with the standard transformation formula for the ADM mass under conformal changes, we obtain
	\[
	\bar m = m + Q = 0.
	\]
	Therefore, equality holds in Herzlich’s theorem, and $(M^n,\bar g)$ is isometric to the complement of a round ball in Euclidean space, say $(\mathbb{R}^n\setminus B_{r_0}(0),\delta)$, where
	\[
	r_0= \frac{(n-1)V_0^{\frac{1}{n-2}}}{\mathcal{H}^-}.
	\] 
	Finally it follows from~\eqref{EquationU-gbar} that the function $U:=V^{-1}$ satisfies the boundary value problem
	\[
	\begin{cases}
		\Delta_\delta U & = 0 \text{ on } \mathbb{R}^n\setminus B_{r_0}(0) \\
		U_{|\partial B_{r_0}(0)} & = V_0^{-1}\\
		U & \rightarrow 1 \text{ as }r\rightarrow\infty.
	\end{cases}
	\]
	Since the maximum principle ensures uniqueness, the solution $U$ must be 
	\[
	U(x)=1+\frac{C}{r^{n-2}}
	\]
	where $C=r_0^{n-2}(V_0^{-1}-1)$. Comparing with the asymptotic expansion of $U$ at infinity, we obtain $C=m$, and hence
	\[
	U(x)=U_m(x),\qquad
	g=U^{\frac{2}{n-2}}_m\delta=g_m,
	\qquad
	E^\flat=\frac{1}{n-2}d\ln U_m=E^\flat_m.
	\]
	Therefore, $(M^n,g,E)$ is isometric to the exterior region $\Omega_{m,r_0}$ of an extremal Reissner--Nordstr\"om manifold of mass $m$ as defined in Section~\ref{MP-Example}.
	
	It remains to rule out the presence of additional asymptotically flat ends. In this situation, the Witten spinor $\Phi\in\Gamma(S_g)$ associated with a chosen end is constructed by prescribing a nonzero constant spinor there and zero boundary conditions at infinity on the remaining ends. In the equality case, this spinor becomes a charged parallel spinor. Hence it is asymptotic to a nonzero constant spinor on the distinguished end and tends to zero on all the others. However, using the asymptotic decay of the charged connection, one shows that a charged parallel spinor cannot decay to zero at infinity unless it vanishes identically in these ends. Since a nontrivial charged parallel spinor has no zero, this yields a contradiction. Consequently, the equality case forces the distinguished end to be the unique asymptotically flat end.
	\qed
	
	\begin{remark}
		The inequality $m\geq |Q|$ holds under the more general assumptions of Proposition~\ref{ExistenceBVP}. However, in the equality case, up to replacing $E$ by $-E$, one can only conclude that $(M^n,g,E)$ admits a negative charged parallel spinor, that the boundary is totally umbilic, and that $\mathcal{H}^-$ is constant along $\partial M$.
	\end{remark}
	
	\begin{remark}
		Instead of using the B\"ar--Hijazi inequality to control the first eigenvalue of the extrinsic Dirac operator, one can apply the Friedrich inequality~\cite{Friedrich3} which states that 
		\begin{equation}\label{FriedrichInequality}
			\lambda^2_1\geq\frac{n-1}{4(n-2)}\inf_{\partial M} \Ss
		\end{equation}
		with equality if and only if $(\partial M,h)$ carries a real Killing spinor. Assuming that the scalar curvature $\Ss>0$, the assumption~\eqref{YamabeMC} is replaced by 
		\[
		\sup_{\partial M}\left(H+(n-1)|g(E,\nu)|\right)\leq \sqrt{\frac{n-1}{n-2}\inf_{\partial M} \Ss}.
		\]
		Under this alternative boundary assumption, the conclusion of Theorem~\ref{ChargedPMT-Connected boundary} is unchanged. In particular, if $\div(E)=0$, the rigidity statement in the equality case remains the same.
	\end{remark}

	%%%%%%%%%%%%%%%%%%%%%%%%%%%%%%%%%%%%%%%%%%%%%%%%%%%%%%%%%%%%%%%%%%%%%%%%%%%%
	
	\section{Manifolds with asymptotically cylindrical ends}\label{PMT-ACE}
	
	%%%%%%%%%%%%%%%%%%%%%%%%%%%%%%%%%%%%%%%%%%%%%%%%%%%%%%%%%%%%%%%%%%%%%%%%%%%%
	
	We now turn to the second geometric setting in which non-trivial electric charge may occur, namely manifolds with asymptotically cylindrical ends. In this section, we prove Theorem~\ref{ChargedPMT-Cylindrical} and characterize its equality case.
	
	%%%%%%%%%%%%%%%%%%%%%%%%%%%%%%%%%%%%%%%%%%%%%%%%%%%%%%%%%%%%%%%%%%%%%%%%%%%%%
	\subsection{The mass--charge inequality}
	%%%%%%%%%%%%%%%%%%%%%%%%%%%%%%%%%%%%%%%%%%%%%%%%%%%%%%%%%%%%%%%%%%%%%%%%%%%%%
	
	In contrast with the case of manifolds with compact boundary, no boundary condition is required here. We consider the space $C^1_c(M,S_g)$ of compactly supported spinor fields on $M$, equipped with the norm
	\[
	||\varphi||_{\mp}^2 := \int_{M}\Big(|\nabla^\mp\varphi|^2+\mathcal{R}^\mp|\varphi|^2\Big)\,d\mu,
	\]
	and denote by $\mathbb{H}^\mp$ its Hilbert space completion. As in~\cite{BartnikChrusciel,Raulot16}, the decay assumptions imply a weighted Poincar\'e inequality and a Schr\"odinger--Lichnerowicz estimate for $D^\mp$. Therefore, by~\cite[Theorem 8.9]{BartnikChrusciel}, the operator
	\[
	D^\mp : \mathbb{H}^\mp \longrightarrow L^2(S_g)
	\]
	is an isomorphism. Arguing exactly as in Section~\ref{PMT-CB}, one obtains, for any asymptotically constant spinor $\Phi_0$, there exists a spinor
	$\Phi\in\Gamma(S_g)$ such that
	\[
	D^-\Phi=0
	\quad\text{and}\quad
	\Phi-\Phi_0\in\mathbb H^-.
	\]
	It follows that
	\[
	\frac{n-1}{2}\omega_{n-1}(m+Q)
	=
	\int_M \Big(|\nabla^-\Phi|^2+\mathcal R^-|\Phi|^2\Big)\,d\mu
	\ge 0,
	\]
	and hence $m\ge |Q|$.

	%%%%%%%%%%%%%%%%%%%%%%%%%%%%%%%%%%%%%%%%%%%%%%%%%%%%%%%%%%%%%%%%%%%%%%%%%%%%
	
	\subsection{The equality case}
	
	%%%%%%%%%%%%%%%%%%%%%%%%%%%%%%%%%%%%%%%%%%%%%%%%%%%%%%%%%%%%%%%%%%%%%%%%%%%%
	
	Throughout this section, we assume that $(M^n,g,E)$ satisfies the equality case of Theorem~\ref{ChargedPMT-Cylindrical} and that $Q\leq 0$. In this situation, the spinor \(\Phi\) constructed above is a negative charged parallel spinor, that is $\nabla^-\Phi=0$. By Proposition~\ref{PropCPS-1}, the spinor \(\Phi\) has no zero and its squared norm $V:=|\Phi|^2$ satisfies~\eqref{ElectricFieldSpinor}. In particular \(dE^\flat=0\). From the asymptotically cylindrical assumptions~\eqref{ACDecay}, the electric field can be written as
	\begin{equation}\label{ElectricFielAsymptotic}
		E=E_\infty+O_1(e^{-\alpha t}),\quad\text{with }E_\infty=Y+q\partial_t
	\end{equation}
	where $Y\in\Gamma(T\Sigma)$ and $q\in C^\infty(\Sigma)$, with \((\Sigma^{n-1},h)\) the limiting cross-section of the cylindrical end.
	\begin{lemma}\label{Lemma-q-positive}
		The function \(q\) is a positive constant on \(\Sigma\).
	\end{lemma}
	
	\pf
	Using~\eqref{ElectricFielAsymptotic} and the identity \(dE^\flat=0\), we obtain
	\[
	d_\Sigma Y^\flat=0,
	\qquad
	d_\Sigma q=0
	\]
	where $d_\Sigma$ denotes the exterior derivative on $\Sigma$. Hence \(q\) is constant on \(\Sigma\).
	
	It remains to determine the sign of \(q\). Since \(\div(E)=0\), applying the divergence theorem to the domain bounded by a large coordinate sphere \(S_r\) in the asymptotically flat end and by a slice \(\Sigma_t\) in the cylindrical end gives
	\[
	0=
	\int_{S_r}g(E,\nu_r)\,d\sigma_r
	-
	\int_{\Sigma_t}g(E,\nu_t)\,d\sigma_t,
	\]
	where \(\nu_t\) denotes the unit normal to \(\Sigma_t\) pointing toward the asymptotically flat end. Letting first \(r\to+\infty\) and then \(t\to+\infty\), and using the definition of the total charge together with the cylindrical asymptotics, we get
	\[
	0=\omega_{n-1}Q+q\,{\rm Vol}(\Sigma,h).
	\]
	Since \(Q\leq0\), we deduce that \(q\geq0\). 
	
	Assume now by contradiction that \(q=0\). Then the preceding identity gives \(Q=0\). Since we are in the equality case of the mass--charge inequality, \(m=|Q|=0\). Moreover, by Proposition~\ref{PropCPS-1} and the assumption \(\div(E)=0\), we have
	\[
	R=(n-1)(n-2)|E|^2\ge0.
	\]
	Thus \((M^n,g)\) is a complete spin manifold with nonnegative scalar curvature and with an asymptotically flat end of zero ADM mass. By the rigidity statement in the positive mass theorem with arbitrary ends due to Cecchini and Zeidler~\cite{CecchiniZeidler}, \((M^n,g)\) is isometric to Euclidean space. This contradicts the existence of an asymptotically cylindrical end. Hence \(q>0\).
	\qed
	
	We now investigate the geometric consequences of the existence of the charged parallel spinor along the cylindrical end. The key point is to extract a non-trivial limiting spinor on the cross-section \((\Sigma,h)\).
	
	In order to analyze the asymptotic behaviour of \(\Phi\), we first identify the spinor bundles over the slices \((\Sigma^{n-1},g_t)\), where \(g_t:=g_{|\Sigma_t}\), with the fixed extrinsic spinor bundle over \((\Sigma^{n-1},h)\). Let
	\[
	B_t:T\Sigma\longrightarrow T\Sigma
	\]
	be the positive \(h\)-symmetric automorphism characterized by
	\[
	g_t(X,Y)=h(B_tX,B_tY).
	\]
	Since \(g_t\to h\) in \(C^2\) by the asymptotically cylindrical assumptions~\eqref{ACDecay}, we have \(B_t\to\Id\) in \(C^2\). By the Bourguignon--Gauduchon construction~\cite{BourguignonGauduchon} recalled in Section~\ref{Spinors}, the maps \(B_t\) induce fiberwise unitary bundle isomorphisms
	\[
	\mathcal I_t:\RSB_{g_t}\longrightarrow \RSB_h,
	\]
	where \(\RSB_h\) denotes the extrinsic spinor bundle associated with \((\Sigma^{n-1},h)\) as defined in~\eqref{SpinorBundleRestriction}. Denote by \(\mult_t\) and \(\nb^{\,t}\) the extrinsic Clifford multiplication and the extrinsic spin connection associated with \(g_t\), and by \(\mult_h\) and \(\nb^h\) the corresponding Clifford multiplication and spin connection on \((\Sigma,h)\), defined in~\eqref{CliffordLCI-Identifications}. Then the identifications \(\mathcal I_t\) preserve Clifford multiplication in the sense that
	\[
	\mathcal I_t\bigl(\mult_t(X)\psi\bigr)
	=
	\mult_h(B_t^{-1}X)\mathcal I_t(\psi)
	\]
	for $X\in\Gamma(T\Sigma)$ and $\psi\in\Gamma(\RSB_{g_t})$. Moreover, transporting the extrinsic spin connections to \(\RSB_h\), we obtain
	\begin{equation}\label{IdentificationLCC}
		\mathcal I_t\circ\nb^{\,t}\circ\mathcal I_t^{-1}
		=
		\nb^h+\mathcal B_t,
	\end{equation}
	where \(\mathcal B_t\in\Omega^1(\Sigma,\End(\RSB_h))\) is locally given, with respect to any local \(h\)-orthonormal frame \((e_i)\), by
	\[
	\mathcal B_t(X)\psi
	=
	\frac14
	\sum_{i,j=1}^{n-1}
	h\Big(
	\big(\widetilde\nabla^{\,t}_X-\nabla^h_X\big)e_i,
	e_j
	\Big)
	\mult_h(e_i)\mult_h(e_j)\psi,
	\]
	with
	\[
	\widetilde\nabla^{\,t}_X Y
	:=
	B_t\bigl(\nabla^{g_t}_X(B_t^{-1}Y)\bigr).
	\]
	Since \(g_t=h+O_2(e^{-\alpha t})\), we have \(\widetilde\nabla^{\,t}-\nabla^h=O_1(e^{-\alpha t})\), and therefore
	\[
	|\mathcal B_t|_h+|\nabla^h\mathcal B_t|_h
	=
	O(e^{-\alpha t}).
	\]
	The following lemma will play a crucial role in the proof.
	\begin{lemma}\label{Lemma-limiting-spinor}
		After passing to a subsequence \(t_j\to+\infty\), the one-parameter family of spinors
		\[
		\psi_t
		:=
		\mathcal I_t\left(\Psi_t\right),\qquad \Psi_t:=e^{\frac{n-2}{2}qt}\Phi|_{\Sigma_t}
		\]
		converges in \(C^1\) to a nonzero spinor $\psi_\infty\in\Gamma(\RSB_h)$. Moreover, 
		\[
		\nb_X^h\psi_\infty
		=
		-\frac12\mult_h(X)\mult_h(Y)\psi_\infty
		+\frac12q\,\mult_h(X)\psi_\infty
		-\frac{n-1}{2}h(Y,X)\psi_\infty
		\]
		for every $X\in\Gamma(T\Sigma)$. 
	\end{lemma}
	
	\pf
	Evaluating~\eqref{ElectricFieldSpinor} on \(\partial_t\), and using the asymptotics of \(E\) and \(g\) together with the fact that the limiting tangential field \(Y\) is independent of \(t\), we obtain
	\[
	\partial_t\ln V=-(n-2)q+O_1(e^{-\alpha t}).
	\]
	Integrating this relation yields
	\begin{equation}\label{V-cyl-expansion}
		V(t,y)=e^{-(n-2)qt}v(y)\bigl(1+O_2(e^{-\alpha t})\bigr),
	\end{equation}
	where \(v\) is a smooth positive function on \(\Sigma\). The \(O_2\)-estimate follows from the asymptotic condition~\eqref{ElectricFielAsymptotic}.
	
	We define the renormalized spinors
	\[
	\Psi_t
	:=
	e^{\frac{n-2}{2}qt}\Phi|_{\Sigma_t}.
	\]
	Then
	\[
	|\Psi_t(y)|^2
	=
	v(y)\bigl(1+O_2(e^{-\alpha t})\bigr),
	\]
	so that \((\Psi_t)\) is uniformly bounded from above and below. We now let
	\[
	\psi_t:=\mathcal I_t(\Psi_t)\in\Gamma(\RSB_h).
	\]
	Since \(\mathcal I_t\) is fiberwise unitary, the family \((\psi_t)\) is uniformly bounded from above and below. We now derive the equation satisfied by \(\psi_t\). We choose the unit normal \(\nu_t\) to \(\Sigma_t\) pointing toward the asymptotically flat end. Because of~\eqref{ElectricFielAsymptotic}, this choice gives
	\[
	E=E_t^\top-q\nu_t+O_1(e^{-\alpha t}),
	\qquad
	E_t^\top\to Y\qquad\text{and}\qquad E_t^\perp:=g(E,\nu_t)\to-q.
	\]
	Using the definition~\eqref{CliffordLCI-Extrinsic} of the extrinsic Levi-Civita connection together with~\eqref{NegativeCPS} and~\eqref{IdentificationLCC}, we obtain
	\begin{equation}\label{psi-t-equation}
		\nb_X^h\psi_t
		=
		\mathcal E_t(X)\psi_t,
	\end{equation}
	for every $X\in\Gamma(T\Sigma)$ and where \(\mathcal E_t\in\Omega^1(\Sigma,\End(\RSB_h))\) is given by
	\[
	\begin{aligned}
		\mathcal E_t(X)
		=
		&-\frac12\mult_h(B_t^{-1}X)\mult_h(B_t^{-1}E_t^\top)
		-\frac12 E_t^\perp\,\mult_h(B_t^{-1}X)\\
		&-\frac{n-1}{2}g_t(E_t^\top,X)\Id
		+\frac12\mult_h(B_t^{-1}A_tX)
		-\mathcal B_t(X).
	\end{aligned}
	\]
	Here \(A_t\) denotes the shape operator of \(\Sigma_t\). By the asymptotically cylindrical assumptions,
	\[
	B_t=\Id+O_2(e^{-\alpha t}),\qquad
	E_t^\top=Y+O_1(e^{-\alpha t}),\qquad
	E_t^\perp=-q+O_1(e^{-\alpha t}),
	\]
	and
	\[
	A_t=O_1(e^{-\alpha t}),\qquad
	\mathcal B_t=O_1(e^{-\alpha t}).
	\]
	Therefore \(\mathcal E_t\longrightarrow \mathcal E_\infty\) in \(C^1\), where
	\[
	\mathcal E_\infty(X)
	=
	-\frac12\mult_h(X)\mult_h(Y)
	+\frac12 q\,\mult_h(X)
	-\frac{n-1}{2}h(Y,X)\Id .
	\]
	In particular, the family \((\mathcal E_t)\) is uniformly bounded in \(C^1\).
	
	Since \((\psi_t)\) is uniformly bounded in \(C^0\) and \((\mathcal E_t)\) is uniformly bounded in \(C^1\), equation~\eqref{psi-t-equation} yields a uniform \(C^2\)-bound for \((\psi_t)\). By the Arzelà--Ascoli theorem, after passing to a subsequence \(t_j\to+\infty\), \(\psi_{t_j}\longrightarrow\psi_\infty\) in \(C^1\). The uniform lower bound on \(|\psi_t|\) implies that the limiting spinor \(\psi_\infty\) is nonzero. Passing to the limit in~\eqref{psi-t-equation}, we obtain
	\[
	\nb_X^h\psi_\infty
	=
	-\frac12\mult_h(X)\mult_h(Y)\psi_\infty
	+\frac12q\,\mult_h(X)\psi_\infty
	-\frac{n-1}{2}h(Y,X)\psi_\infty
	\]
	for $X\in\Gamma(T\Sigma)$. 
	\qed 
	
	We now show that the tangential component of the limiting electric field must vanish identically.
	\begin{lemma}\label{Lemma-Y-zero}
		The tangential component of the limiting electric field vanishes, namely \(Y\equiv0\). Consequently, the limiting spinor satisfies
		\[
		\nb_X^h\psi_\infty
		=
		\frac 12q\,\mult_h(X)\psi_\infty,
		\]
		so that \(\psi_\infty\) is an extrinsic Killing spinor on the limiting cross-section \((\Sigma^{n-1},h)\).
	\end{lemma}
	
	\pf
	The extrinsic Dirac operator on \(\RSB_h\) is defined by~\eqref{ExtrinsicDirac}. Using Lemma~\ref{Lemma-limiting-spinor}, we compute
	\[
	\begin{aligned}
		\D_h\psi_\infty
		&=
		-\frac12
		\sum_{i=1}^{n-1}
		\mult_h(e_i)\mult_h(e_i)\mult_h(Y)\psi_\infty
		+
		\frac q2
		\sum_{i=1}^{n-1}
		\mult_h(e_i)\mult_h(e_i)\psi_\infty\\
		&\quad
		-\frac{n-1}{2}
		\sum_{i=1}^{n-1}
		h(Y,e_i)\mult_h(e_i)\psi_\infty
	\end{aligned}
	\]
	so that
	\[
	\D_h\psi_\infty
	=
	-\frac{n-1}{2}q\,\psi_\infty.
	\]
	Thus \(\psi_\infty\) is a nonzero eigenspinor of \(\D_h\) with eigenvalue $-(n-1)q/2$. On the other hand, since 
	\[
	R=(n-1)(n-2)|E|^2,
	\] 
	the Gauss equation on the slices \(\Sigma_t\) yields
	\[
	\Ss_t
	=
	(n-1)(n-2)|E|^2
	-
	2{\rm Ric}(\nu_t,\nu_t)
	+
	H_t^2-|A_t|^2
	\]
	where $\Ss_t$ denotes the scalar curvature of \((\Sigma^{n-1},g_t)\). Moreover, the asymptotically cylindrical assumptions imply
	\[
	{\rm Ric}(\nu_t,\nu_t)=O(e^{-\alpha t}).
	\]
	Since \(g_t\to h\) in \(C^2\), we have \(\Ss_t\rightarrow \Ss_h\) where $\Ss_h$ denotes the scalar curvature of the limiting cross-section \((\Sigma^{n-1},h)\). Furthermore, by~\eqref{ElectricFielAsymptotic},
	\[
	E=Y+q\partial_t+O_1(e^{-\alpha t}),
	\]
	and therefore
	\[
	|E|^2
	=
	|Y|_h^2+q^2+O_1(e^{-\alpha t}).
	\]
	Passing to the limit as \(t\to+\infty\) in the Gauss equation yields
	\begin{equation}\label{Rh-cylinder}
		\Ss_h=(n-1)(n-2)(|Y|_h^2+q^2).
	\end{equation}
	Now recall that the spectrum of the extrinsic Dirac operator is obtained by symmetrizing the spectrum of the intrinsic Dirac operator on \((\Sigma^{n-1},h)\) (see for instance~\cite{Ginoux}). Friedrich's inequality~\eqref{FriedrichInequality} on the limiting cross-section \((\Sigma^{n-1},h)\) therefore gives
	\[
	\left(\frac{n-1}{2}q\right)^2
	\ge
	\frac{n-1}{4(n-2)}\inf_\Sigma \Ss_h.
	\]
	Using~\eqref{Rh-cylinder}, this yields
	\[
	\inf_\Sigma |Y|_h^2=0.
	\]
	Moreover equality holds in Friedrich's inequality and in particular \(\Ss_h\) is constant. Since \(q\) is constant, it follows from~\eqref{Rh-cylinder} that \(|Y|_h^2\) is constant. Together with \(\inf_\Sigma |Y|_h^2=0\), this gives $Y\equiv0$. Substituting this into the equation satisfied by \(\psi_\infty\) in Lemma~\ref{Lemma-limiting-spinor}, we conclude that \(\psi_\infty\) is an extrinsic Killing spinor on \((\Sigma^{n-1},h)\) with Killing constant \(q/2\).
	\qed
	
	We now study the conformal metric
	\[
	\overline g:=V^{\frac{2}{n-2}}g,
	\]
	which will allow us to complete the rigidity argument.
	\begin{lemma}\label{Lemma-conformal-metric}
		The metric \(\bar g\) is asymptotically flat on the distinguished end, Ricci-flat, and has ADM mass $\bar m=0$. Moreover, there exist \(A>0\) and a coordinate chart on the asymptotically cylindrical end with values in $(0,A]\times\Sigma$ in which
		\[
		\bar g=g_C+\kappa,
		\]
		where
		\[
		g_C=d\rho^2+\rho^2(q^2h)
		\]
		and
		\[
		|(\nabla^C)^j\kappa|_{g_C}=O(\rho^{\beta-j}),
		\qquad j=0,1,2,
		\]
		for some \(\beta>0\). Here \(\nabla^C\) denotes the Levi--Civita connection of \(g_C\).
	\end{lemma}
	
	\pf
	By Proposition~\ref{Conformal-PCS}, the metric \(\bar g\) carries a nonzero parallel spinor and is therefore Ricci-flat. Moreover, the same argument as in the proof of Theorem~\ref{ChargedPMT-Connected boundary} (see Section~\ref{ProofTHMCB}) yields
	\[
	V=1+\frac{Q}{r^{n-2}}+o_2(r^{2-n}),
	\]
	so that \((M^n,\bar g)\) is asymptotically flat. The usual transformation law for the ADM mass under conformal changes gives \(\bar m=m+Q\). Since equality means \(m=|Q|=-Q\), we conclude that \(\bar m=0\).
	
	To describe the geometry of \(\bar g\) near the cylindrical end, we differentiate~\eqref{V-cyl-expansion} tangentially along \(\Sigma\) and use~\eqref{ElectricFieldSpinor} to obtain
	\[
	d_\Sigma(\ln v)=-(n-2)Y^\flat.
	\]
	By Lemma~\ref{Lemma-Y-zero}, \(Y\equiv0\), and hence \(v\) is constant. Therefore
	\[
	V(t,y)=C e^{-(n-2)qt}\bigl(1+O_2(e^{-\alpha t})\bigr)
	\]
	for some constant \(C>0\). Since \(q>0\) by Lemma~\ref{Lemma-q-positive}, we may introduce the coordinate
	\[
	\rho:=\frac{C^{1/(n-2)}}{q}e^{-qt}.
	\]
	Then \(\rho\to0\) as \(t\to+\infty\) and
	\[
	e^{-\alpha t}=O(\rho^\beta),
	\qquad
	\beta:=\frac{\alpha}{q}>0.
	\]
	Consequently,
	\[
	V^{\frac{2}{n-2}}
	=
	q^2\rho^2\bigl(1+O_2(\rho^\beta)\bigr).
	\]
	Combining this with
	\[
	g=dt^2+h+O_2(e^{-\alpha t})
	=
	dt^2+h+O_2(\rho^\beta)
	\]
	and the identity \(d\rho=-q\rho\,dt\), we obtain
	\[
	\bar g
	=
	d\rho^2+\rho^2(q^2h)+\kappa,
	\]
	where \(\kappa\) satisfies the required estimates.
	\qed
	
	We now complete the proof of the rigidity statement. Let \(\widehat M\) be the space obtained from \(M\) by adding one point \(p\) corresponding to the asymptotically cylindrical end. By Lemma~\ref{Lemma-conformal-metric}, \(\bar g\) extends to \(\widehat M\) as an asymptotically flat metric with an isolated conical singularity at \(p\) in the sense of Dai--Sun--Wang~\cite{DaiSunWang1,DaiSunWang2}. Moreover, \(\bar g\) is Ricci-flat on \(M\) and has zero ADM mass.
	
	Consequently, \((\widehat M,\bar g)\) satisfies the assumptions of the positive mass theorem with isolated conical singularities of~\cite{DaiSunWang1,DaiSunWang2}. By the rigidity statement, \((\widehat M,\bar g)\) is isometric to Euclidean space. Choosing the isometry so that \(p\) corresponds to the origin, we identify 
	\[
	(M^n,\bar g)\simeq(\mathbb R^n\setminus\{0\},\delta).
	\]
	
	Let \(U:=V^{-1}\). Then
	\[
	g=U^{\frac{2}{n-2}}\delta,
	\qquad
	E^\flat=\frac1{n-2}d\ln U.
	\]
	Moreover, by~\eqref{EquationU-gbar} and the assumption \(\div(E)=0\),
	\[
	\Delta_\delta U=0
	\qquad\text{on }\mathbb R^n\setminus\{0\}.
	\]  
	On the asymptotically flat end,
	\[
	U
	=
	1-\frac{Q}{r^{n-2}}+o_2(r^{2-n})
	=
	1+\frac{m}{r^{n-2}}+o_2(r^{2-n}),
	\]
	since \(m=|Q|=-Q>0\). The asymptotic behaviour of \(V\) at the cylindrical end implies that \(U\) has an isolated pole at \(0\). Hence \(U\) is a positive harmonic function on \(\mathbb R^n\setminus\{0\}\), tends to \(1\) at infinity, and has an isolated pole at \(0\). By Bôcher's theorem,
	\[
	U(x)=1+\frac{a}{|x|^{n-2}}
	\]
	for some \(a>0\). Comparing with the asymptotic expansion at infinity yields \(a=m\), and therefore \(U=U_m\). Consequently,
	\[
	g=U_m^{\frac{2}{n-2}}\delta=g_m,
	\qquad
	E^\flat=\frac1{n-2}d\ln U_m=E_m^\flat.
	\]
	Therefore \((M^n,g,E)\) is isometric to the extremal Reissner--Nordstr\"om manifold of mass \(m>0\) and negative charge.

	%%%%%%%%%%%%%%%%%%%%%%%%%%%%%%%%%%%%%%%%%%%%%%%%%%%%%%%%%%%%%%%%%%%%%%%%%%%%%%%%%%%%

	\bibliographystyle{alpha}     
	\bibliography{BiblioHabilitation}

@article {Bar2,
    AUTHOR = {Bär, C.},
     TITLE = {Real {K}illing spinors and holonomy},
   JOURNAL = {Comm. Math. Phys.},
  FJOURNAL = {Communications in Mathematical Physics},
    VOLUME = {154},
      YEAR = {1993},
    NUMBER = {3},
     PAGES = {509--521},
}

@article {Bar3,
    AUTHOR = {Bär, C.},
     TITLE = {Lower eigenvalue estimates for {D}irac operators},
   JOURNAL = {Math. Ann.},
  FJOURNAL = {Mathematische Annalen},
    VOLUME = {293},
      YEAR = {1992},
    NUMBER = {1},
     PAGES = {39--46},
}

@article {Bartnik1,
    AUTHOR = {Bartnik, R.},
     TITLE = {The mass of an asymptotically flat manifold},
   JOURNAL = {Comm. Pure Appl. Math.},
  FJOURNAL = {Communications on Pure and Applied Mathematics},
    VOLUME = {39},
      YEAR = {1986},
    NUMBER = {5},
     PAGES = {661--693},
}

@article {BartnikChrusciel,
    AUTHOR = {Bartnik, R. and Chru\'{s}ciel, P. T.},
     TITLE = {Boundary value problems for {D}irac-type equations},
   JOURNAL = {J. Reine Angew. Math.},
  FJOURNAL = {Journal f\"{u}r die Reine und Angewandte Mathematik. [Crelle's
              Journal]},
    VOLUME = {579},
      YEAR = {2005},
     PAGES = {13--73},
}

@article {BourguignonGauduchon,
    AUTHOR = {Bourguignon, J.-P. and Gauduchon, P.},
     TITLE = {Spineurs, op\'{e}rateurs de {D}irac et variations de m\'{e}triques},
   JOURNAL = {Comm. Math. Phys.},
  FJOURNAL = {Communications in Mathematical Physics},
    VOLUME = {144},
      YEAR = {1992},
    NUMBER = {3},
     PAGES = {581--599},
}

@book {BourguignonHijaziMilhoratMoroianu,
    AUTHOR = {Bourguignon, J.- P. and Hijazi, O. and Milhorat,
              J.- L. and Moroianu, A. and Moroianu, S.},
     TITLE = {A spinorial approach to {R}iemannian and conformal geometry},
    SERIES = {EMS Monographs in Mathematics},
 PUBLISHER = {European Mathematical Society (EMS), Z\"{u}rich},
      YEAR = {2015},
     PAGES = {ix+452},
      ISBN = {978-3-03719-136-1},
}

@incollection {BrayHirschKazarasKhuri,
    AUTHOR = {Bray, H. L. and Hirsch, S. and Kazaras, D. and Khuri,
              M. and Zhang, Y.},
     TITLE = {Spacetime harmonic functions and applications to mass},
 BOOKTITLE = {Perspectives in scalar curvature. {V}ol. 2},
     PAGES = {593--639},
 PUBLISHER = {World Sci. Publ., Hackensack, NJ},
      YEAR = { 2023},
      ISBN = {978-981-124-999-0; 978-981-124-935-8; 978-981-124-936-5},
   MRCLASS = {53C21 (53C20 53C43)},
  MRNUMBER = {4577926},
}

@incollection {Chrusciel1,
    AUTHOR = {Chru\'{s}ciel, P. T.},
     TITLE = {Boundary conditions at spatial infinity from a {H}amiltonian
              point of view},
 BOOKTITLE = {Topological properties and global structure of space-time
              ({E}rice, 1985)},
    SERIES = {NATO Adv. Sci. Inst. Ser. B Phys.},
    VOLUME = {138},
     PAGES = {49--59},
 PUBLISHER = {Plenum, New York},
      YEAR = {1986},
}

@article {ChruscielRealTod,
    AUTHOR = {Chru\'sciel, P. T. and Reall, H. S. and Tod, P.},
     TITLE = {On {I}srael-{W}ilson-{P}erj\'es black holes},
   JOURNAL = {Classical Quantum Gravity},
  FJOURNAL = {Classical and Quantum Gravity},
    VOLUME = {23},
      YEAR = {2006},
    NUMBER = {7},
     PAGES = {2519--2540},
      ISSN = {0264-9381,1361-6382},
   MRCLASS = {83C57},
  MRNUMBER = {2215078},
MRREVIEWER = {Piero\ Nicolini},
       DOI = {10.1088/0264-9381/23/7/018},
       URL = {https://doi.org/10.1088/0264-9381/23/7/018},
}

@article{CecchiniZeidler,
 author = {Cecchini, S. and Zeidler, R.},
 title = {The positive mass theorem and distance estimates in the spin setting},
 fjournal = {Transactions of the American Mathematical Society},
 journal = {Trans. Am. Math. Soc.},
 issn = {0002-9947},
 volume = {377},
 number = {8},
 pages = {5271--5288},
 year = {2024},
}

@article{CruzLimaDeSousa,
 author = {Cruz, T. and Lima, V. and de Sousa, A.},
 title = {Min-max minimal surfaces, horizons and electrostatic systems},
 fjournal = {Journal of Differential Geometry},
 journal = {J. Differ. Geom.},
 issn = {0022-040X},
 volume = {128},
 number = {2},
 pages = {583--637},
 year = {2024},
 language = {English},
 doi = {10.4310/jdg/1727712890},
 url = {projecteuclid.org/journals/journal-of-differential-geometry/volume-128/issue-2/Min-max-minimal-surfaces-horizons-and-electrostatic-systems/10.4310/jdg/1727712890.full},
 zbMATH = {7939776},
 Zbl = {1552.53009}
}

@article {DaiSunWang1,
    AUTHOR = {Dai, X. and Sun, Y. and Wang, C.},
     TITLE = {The positive mass theorem for asymptotically flat manifolds with isolated conical singularities},
   JOURNAL = {Sci. China Math. },
   YEAR = {2024},
  URL = {https://doi.org/10.1016/j.difgeo.2014.09.003},
 }

@article{DaiSunWang2,
 author = {Dai, X. and Sun, Y. and Wang, C.},
 title = {Positive mass theorem for asymptotically flat spin manifolds with isolated conical singularities},
 fjournal = {Transactions of the American Mathematical Society},
 journal = {Trans. Am. Math. Soc.},
 issn = {0002-9947},
 volume = {378},
 number = {4},
 pages = {2617--2642},
 year = {2025},
 language = {English},
 doi = {10.1090/tran/9331},
 zbMATH = {8013995},
 Zbl = {1566.53055}
}

@article{FreitasLeandroRibeiroSabo,
 author = {Freitas, A. and Leandro, B. and Ribeiro, E. and Sabo, G.},
 title = {Geometric inequalities for electrostatic systems with boundary},
 fjournal = {Classical and Quantum Gravity},
 journal = {Classical Quantum Gravity},
 issn = {0264-9381},
 volume = {43},
 number = {3},
 pages = {32},
 note = {Id/No 035003},
 year = {2026},
 language = {English},
 doi = {10.1088/1361-6382/ae3a2f},
 zbMATH = {8168497}
}

@article {Friedrich3,
    AUTHOR = {Friedrich, T.},
     TITLE = {Der erste {E}igenwert des {D}irac-{O}perators einer kompakten,
              {R}iemannschen {M}annigfaltigkeit nichtnegativer
              {S}kalarkr\"{u}mmung},
   JOURNAL = {Math. Nachr.},
  FJOURNAL = {Mathematische Nachrichten},
    VOLUME = {97},
      YEAR = {1980},
     PAGES = {117--146},
}

@article {GibbonsHawkingHorowitzPerry,
    AUTHOR = {Gibbons, G. W. and Hawking, S. W. and Horowitz, G. T. and
              Perry, M. J.},
     TITLE = {Positive mass theorems for black holes},
   JOURNAL = {Comm. Math. Phys.},
  FJOURNAL = {Communications in Mathematical Physics},
    VOLUME = {88},
      YEAR = {1983},
    NUMBER = {3},
     PAGES = {295--308},
}

@article {GibbonsHull,
    AUTHOR = {Gibbons, G. W. and Hull, C. M.},
     TITLE = {A {B}ogomolny bound for general relativity and solitons in
              {$N=2$}\ supergravity},
   JOURNAL = {Phys. Lett. B},
  FJOURNAL = {Physics Letters. B. Particle Physics, Nuclear Physics and
              Cosmology},
    VOLUME = {109},
      YEAR = {1982},
    NUMBER = {3},
     PAGES = {190--194},
      ISSN = {0370-2693,1873-2445},
   MRCLASS = {83C30 (53A99 83E50)},
  MRNUMBER = {644079},
MRREVIEWER = {M.\ J.\ Perry},
       DOI = {10.1016/0370-2693(82)90751-1},
       URL = {https://doi.org/10.1016/0370-2693(82)90751-1},
}

@book {Ginoux,
    AUTHOR = {Ginoux, N.},
     TITLE = {The {D}irac spectrum},
    SERIES = {Lecture Notes in Mathematics},
    VOLUME = {1976},
 PUBLISHER = {Springer-Verlag, Berlin},
      YEAR = {2009},
     PAGES = {xvi+156},
      ISBN = {978-3-642-01569-4},
}

@article {Herzlich1,
    AUTHOR = {Herzlich, M.},
     TITLE = {A {P}enrose-like inequality for the mass of {R}iemannian
              asymptotically flat manifolds},
   JOURNAL = {Comm. Math. Phys.},
  FJOURNAL = {Communications in Mathematical Physics},
    VOLUME = {188},
      YEAR = {1997},
    NUMBER = {1},
     PAGES = {121--133},
}

@incollection {Herzlich2,
    AUTHOR = {Herzlich, M.},
     TITLE = {Minimal surfaces, the {D}irac operator and the {P}enrose
              inequality},
 BOOKTITLE = {S\'{e}minaire de {T}h\'{e}orie {S}pectrale et {G}\'{e}om\'{e}trie, {V}ol. 20,
              {A}nn\'{e}e 2001--2002},
    SERIES = {S\'{e}min. Th\'{e}or. Spectr. G\'{e}om.},
    VOLUME = {20},
     PAGES = {9--16},
 PUBLISHER = {Univ. Grenoble I, Saint-Martin-d'H\`eres},
      YEAR = {2002},
}

@article {Hijazi1,
    AUTHOR = {Hijazi, O.},
     TITLE = {Premi\`ere valeur propre de l'op\'{e}rateur de {D}irac et nombre de
              {Y}amabe},
   JOURNAL = {C. R. Acad. Sci. Paris S\'{e}r. I Math.},
  FJOURNAL = {Comptes Rendus de l'Acad\'{e}mie des Sciences. S\'{e}rie I.
              Math\'{e}matique},
    VOLUME = {313},
      YEAR = {1991},
    NUMBER = {12},
     PAGES = {865--868},
}

@article {Hijazi2,
    AUTHOR = {Hijazi, O.},
     TITLE = {A conformal lower bound for the smallest eigenvalue of the
              {D}irac operator and {K}illing spinors},
   JOURNAL = {Comm. Math. Phys.},
  FJOURNAL = {Communications in Mathematical Physics},
    VOLUME = {104},
      YEAR = {1986},
    NUMBER = {1},
     PAGES = {151--162},
}

@article {HijaziMontiel1,
    AUTHOR = {Hijazi, O. and Montiel, S.},
     TITLE = {Extrinsic {K}illing spinors},
   JOURNAL = {Math. Z.},
  FJOURNAL = {Mathematische Zeitschrift},
    VOLUME = {244},
      YEAR = {2003},
    NUMBER = {2},
     PAGES = {337--347},
}

@article{Innami,
 author = {Innami, N.},
 title = {Splitting theorems of {Riemannian} manifolds},
 fjournal = {Compositio Mathematica},
 journal = {Compos. Math.},
 issn = {0010-437X},
 volume = {47},
 pages = {237--247},
 year = {1982},
 language = {English},
 url = {https://eudml.org/doc/89573},
 zbMATH = {3812455},
 Zbl = {0514.53040}
}

@article {Jaracz,
    AUTHOR = {Jaracz, J. S.},
     TITLE = {The {P}enrose inequality and positive mass theorem with charge
              for manifolds with asymptotically cylindrical ends},
   JOURNAL = {Ann. Henri Poincar\'e},
  FJOURNAL = {Annales Henri Poincar\'e. A Journal of Theoretical and
              Mathematical Physics},
    VOLUME = {21},
      YEAR = {2020},
    NUMBER = {8},
     PAGES = {2581--2609},
      ISSN = {1424-0637,1424-0661},
   MRCLASS = {83C22 (53C20 83C40)},
  MRNUMBER = {4127376},
MRREVIEWER = {Veselin\ T.\ Videv},
       DOI = {10.1007/s00023-020-00927-z},
       URL = {https://doi.org/10.1007/s00023-020-00927-z},
}

@book {LawsonMichelsohn,
    AUTHOR = {Lawson, Jr., H. B. and Michelsohn, M.-L.},
     TITLE = {Spin geometry},
    SERIES = {Princeton Mathematical Series},
    VOLUME = {38},
 PUBLISHER = {Princeton University Press, Princeton, NJ},
      YEAR = {1989},
     PAGES = {xii+427},
      ISBN = {0-691-08542-0},
}

@book {Lee,
    AUTHOR = {Lee, D. A.},
     TITLE = {Geometric relativity},
    SERIES = {Graduate Studies in Mathematics},
    VOLUME = {201},
 PUBLISHER = {American Mathematical Society, Providence, RI},
      YEAR = {2019},
     PAGES = {xii+361},
      ISBN = {978-1-4704-5081-6},
}

@article {McCormick,
    AUTHOR = {McCormick, S},
     TITLE = {A conformal reduction for the {X}-{ADM}  mass},
   JOURNAL = {arXiv:2607.03629},
      YEAR = {2026},
}

@article{Raulot16,
 AUTHOR = {Raulot, S.},
 TITLE = {Positive energy theorems for spin initial data with charge},
JOURNAL = {Classical and Quantum Gravity},
    VOLUME = {43},
      YEAR = {2026},
    NUMBER = {4},
     PAGES = {045011},
}

@article {SchoenYau1,
    AUTHOR = {Schoen, R. M. and Yau, S.-T.},
     TITLE = {On the proof of the positive mass conjecture in general
              relativity},
   JOURNAL = {Comm. Math. Phys.},
  FJOURNAL = {Communications in Mathematical Physics},
    VOLUME = {65},
      YEAR = {1979},
    NUMBER = {1},
     PAGES = {45--76},
}

@article {SchoenYau6,
    AUTHOR = {Schoen, R. M. and Yau, S.-T.},
     TITLE = {The energy and the linear momentum of space-times in general
              relativity},
   JOURNAL = {Comm. Math. Phys.},
  FJOURNAL = {Communications in Mathematical Physics},
    VOLUME = {79},
      YEAR = {1981},
    NUMBER = {1},
     PAGES = {47--51},
}

@article {Tod83,
    AUTHOR = {Tod, K. P.},
     TITLE = {All metrics admitting super-covariantly constant spinors},
   JOURNAL = {Phys. Lett. B},
  FJOURNAL = {Physics Letters. B. Particle Physics, Nuclear Physics and
              Cosmology},
    VOLUME = {121},
      YEAR = {1983},
    NUMBER = {4},
     PAGES = {241--244},
      ISSN = {0370-2693,1873-2445},
   MRCLASS = {83C15 (53C80 83C20)},
  MRNUMBER = {690024},
MRREVIEWER = {Gary\ T.\ Horowitz},
       DOI = {10.1016/0370-2693(83)90797-9},
       URL = {https://doi.org/10.1016/0370-2693(83)90797-9},
}

@article {WangMcK1,
    AUTHOR = {Wang, M. Y.},
     TITLE = {Parallel spinors and parallel forms},
   JOURNAL = {Ann. Global Anal. Geom.},
  FJOURNAL = {Annals of Global Analysis and Geometry},
    VOLUME = {7},
      YEAR = {1989},
    NUMBER = {1},
     PAGES = {59--68},
}

@article {WangMcK2,
    AUTHOR = {Wang, M. Y.},
     TITLE = {On non-simply connected manifolds with non-trivial parallel
              spinors},
   JOURNAL = {Ann. Global Anal. Geom.},
  FJOURNAL = {Annals of Global Analysis and Geometry},
    VOLUME = {13},
      YEAR = {1995},
    NUMBER = {1},
     PAGES = {31--42},
}

@article {Witten1,
    AUTHOR = {Witten, E.},
     TITLE = {A new proof of the positive energy theorem},
   JOURNAL = {Comm. Math. Phys.},
  FJOURNAL = {Communications in Mathematical Physics},
    VOLUME = {80},
      YEAR = {1981},
    NUMBER = {3},
     PAGES = {381--402},
}

	%%%%%%%%%%%%%%%%%%%%%%%%%%%%%%%%%%%%%%%%%%%%%%%%%%%%%%%%%%%%%%%%%%%%%%%%%%%%%%

\end{document}